\documentclass[11pt]{article}

\usepackage{amsmath,amsthm,amssymb}

\usepackage[francais,english]{babel}
\usepackage[T1]{fontenc}
\usepackage{latexsym}
\usepackage{subfigure, epsfig, epsf}
\usepackage{color,fancyhdr,lastpage,graphicx}
\usepackage{xspace}
\usepackage{setspace}
\usepackage{bbm}

\newcommand{\NN}{\mathbb{N}}    \newcommand{\QQ}{\mathbb{Q}}  \newcommand{\HH}{\mathbb{H}}
\newcommand{\RR}{\mathbb{R}}    \newcommand{\LL}{\mathbb{L}}  \newcommand{\WW}{\mathbb{W}}  
\newcommand{\MMM}{\mathbb{M}} \newcommand{\OO}{\mathbb{O}}  \newcommand{\VV}{\mathbb{V}}

 \newcommand{\mcC}{\mathcal{C}} \newcommand{\mcF}{\mathcal{F}} \newcommand{\mcU}{\mathcal{U}}
\newcommand{\mcH}{\mathcal{H}} \newcommand{\mcalV}{\mathcal{V}} \newcommand{\mcA}{\mathcal{A}} \newcommand{\mcB}{\mathcal{B}}

\newcommand{\ff}{{\frak{z}}}

\newcommand{\bfg}{{\bf g}}   \newcommand{\bfi}{{\bf i}}   
  \newcommand{\be}{{\bf e}} \newcommand{\ba}{{\bf a}}

\newcommand{\bfc}{{\bf c}}

 \newcommand{\la}{\leftarrow}

\newcommand{\EE}{\mathbb{E}} \newcommand{\PP}{\mathbb{P}}

\newcommand{\bm}{Brownian motion }

\newcommand{\rv}{random variable }
    \newcommand{\rvs}{random variables }

\newcommand{\as}{as }

    \newcommand{\sdes}{stochastic differential equations }
\newcommand{\rw}{random walk }
    \newcommand{\rws}{random walks }

\newcommand{\Vol}{\small{\textsc{Vol}}}

\newcommand{\LEB}{\textsc{Leb}}

\newcommand{\st}{such that }

\newcommand{\varep}{\varepsilon}

\renewcommand{\leq}{\leqslant}
\renewcommand{\geq}{\geqslant}

\renewcommand{\la}{\lambda}
\newcommand{\al}{\alpha}

\newcommand{\wrt}{with respect to }
\renewcommand{\st}{such that }
\newcommand{\bel}{\overset{\leftarrow}{\be}}
\newcommand{\wlg}{without loss of generality }

\newtheorem{thm}{ {\sc Theorem} }

\newtheorem{conj}[thm]{ {\sc Conjecture} }

\newtheorem{lem}[thm]{ {\sc Lemma} }

\newtheorem{prop}[thm]{ {\sc Proposition}}

\newtheorem{defn}[thm]{ {\sc Definition}}

\newtheorem{PropDefn}[thm]{ {\sc Proposition/Definition} }
\newtheorem{ThmDefn}[thm]{ {\sc Theorem/Definition} }

\newtheorem{OpenPb}[thm]{ {\sc Open Problem} }

\numberwithin{equation}{section} 

\newenvironment{Dem}{%
    \begin{list}{$\lhd$ {\sc Proof --}}{%
        \setlength{\topsep}{0pt}%
        \setlength{\leftmargin}{0pt}%
        \setlength{\rightmargin}{0pt}%
        \setlength{\listparindent}{0pt}%
        \setlength{\itemindent}{0pt}%
        \setlength{\parsep}{0pt}%
        \addtolength{\leftmargin}{20pt}%
        \addtolength{\rightmargin}{0pt}%
    } \item }{\hfill{\space $\rhd$}\end{list}\smallskip}

\newenvironment{DemOuv}{%
    \begin{list}{$\lhd$ {\sc Proof --}}{%
        \setlength{\topsep}{0pt}%
        \setlength{\leftmargin}{0pt}%
        \setlength{\rightmargin}{0pt}%
        \setlength{\listparindent}{0pt}%
        \setlength{\itemindent}{0pt}%
        \setlength{\parsep}{0pt}%
        \addtolength{\leftmargin}{20pt}%
        \addtolength{\rightmargin}{0pt}%
    } \item }{\hfill{\space }\end{list}\smallskip}

\newenvironment{DemFerm}{%
    \begin{list}{ }{%
        \setlength{\topsep}{0pt}%
        \setlength{\leftmargin}{0pt}%
        \setlength{\rightmargin}{0pt}%
        \setlength{\listparindent}{0pt}%
        \setlength{\itemindent}{0pt}%
        \setlength{\parsep}{0pt}%
        \addtolength{\leftmargin}{20pt}%
        \addtolength{\rightmargin}{0pt}%
    } \item }{\hfill{\space $\rhd$}\end{list}\smallskip}

\newenvironment{SousDem}{%
    \begin{list}{$\circ$}{%
        \setlength{\topsep}{0pt}%
        \setlength{\leftmargin}{0pt}%
        \setlength{\rightmargin}{0pt}%
        \setlength{\listparindent}{0pt}%
        \setlength{\itemindent}{0pt}%
        \setlength{\parsep}{0pt}%
        \addtolength{\leftmargin}{15pt}%
        \addtolength{\rightmargin}{0pt}%
    } \item }{\hfill{\space $\circ$}\end{list}\smallskip}

\newcommand*\tfstyle{\bfseries}
\newcommand*\tfafter{\hspace*{0.5ex}}
\newcommand*\tfbefore{}



\title{\LARGE {\sc A pathwise approach to relativistic diffusions}}
\author{Isma\"el Bailleul\footnote{Statistical Laboratory, Center for Mathematical Sciences, Wilberforce Road, Cambridge, CB3 0WB, UK}}

\begin{document}

\maketitle

\selectlanguage{english}

\begin{center}
\Large{Submitted article}
\end{center}

\begin{abstract}
A new class of relativistic diffusions encompassing all the previously studied examples has recently been introduced in the article \cite{ROUPUnifying2} of C. Chevalier and F. Debbasch, both in a heuristic and analytic way. A pathwise approach of these processes is proposed here, in the general framework of Lorentzian geometry. In considering the dynamics of the random motion in strongly causal spacetimes, we are able to give a simple definition of the one-particle distribution function associated with each process of the class and prove its fundamental property. This result not only provides a dynamical justification of the analytical approach developped up to now (enabling us to recover many of the results obtained so far), but it provides a new general $H$-theorem. It also sheds some light on the importance of the large scale structure of the manifold in the asymptotic behaviour of the Franchi-Le Jan process. This pathwise approach is also the source of many interesting questions that have no analytical counterparts.
\end{abstract}

\medskip


\noindent $1991$ \textit{Mathematics Subject Classification}. Primary 60H10, Secondary 83C99.
\medskip

\begin{center}
\tableofcontents
\end{center}

\section{Introduction}
\label{poiuyt}

The present article is at the confluence of two different stories that have met recently.

\medskip

The first was initiated by Dudley in a paper \cite{Dudley1}, written in $1966$, where he describes the class of random Markov timelike paths in Minkowski spacetime whose laws are defined independently of any rest frame. These random paths represent the trajectories of particles whose speed is less than the speed of light, and whose laws are invariant by the action of the isometry group of the space. He proves in this article that there exists essentially a unique way of constructing $\mcC^1$ random paths having the above properties. The phase space $\bigl(\RR\times\RR^3\bigr)\times\HH$ is well adapted to describe it. We write here $\HH$ for the half-unit sphere $\bigl\{\zeta=(t,x)\in\RR\times\RR^3\,;\, q(\zeta):= t^2-|x|_{\textrm{Eucl}}^2 = 1,\,t>0\bigr\}$ of the spacetime $\RR\times\RR^3$, equipped with the quadratic form $q$. The restriction of $q$ to any tangent hyperplane of $\HH$ is definite-negative. Any $\mcC^1$ timelike path whose $t$-co-ordinate increases can always be re-parametrized in such a way that its speed belongs to $\HH$. Random $\mcC^1$ timelike paths $\{\gamma_s\}_{s\geq 0} = \bigr\{\gamma_0+\int_0^s \dot\gamma_r\,dr\bigr\}$ are determined by their $\HH$-valued speed process $\{\dot\gamma_s\}_{s\geq 0}$ which has no other choice than being a \bm on $\HH$ (up to a constant time scaling). Minkowski spacetime thus has a canonical diffusion, in the same way as \bm is canonically associated to Euclidean space. 

This fundamental work had to wait for the development of stochastic analysis and the article \cite{FranchiLeJan} of Franchi and Le Jan, in $2005$, to see its scope extended to the realm of general relativity. They defined a diffusion in any Lorentzian manifold using a stochastic development procedure similar in spirit to the construction of \bm promoted by Malliavin and Elworthy, using \sdes in the orthonormal frame bundle of the manifold.

\medskip

The other story was born immediately after Einstein's theory of relativity and gravitation was accepted and spread in the scientific community. It deals with the extension of Boltzmann theory of gases to the relativistic framework. Although Boltzmann model is primarily a particle model of gases, most of the works have been on understanding the macroscopic behaviour of relativistic gases through the study of the raltivistic Boltzmann equation. One had to wait the nineties and the article \cite{ROUP0} of F. Debbasch, K. Mallick and J.P. Rivet to see the introduction of a probabilistic mesoscopic model of diffusion of a particle in a fluid, under the form of a special relativistic counterpart of Ornstein-Uhlenbeck process. Generalisations of this model to the framework of general relativity have been given in later articles.

\medskip

These two stories have recently met with the proposition, made in the article \cite{ROUPUnifying2} of C. Chevalier and F. Debbasch, to define a class of random processes including Dudley's process and the relativistic Ornstein-Uhlenbeck process, and characterized by the following property. There exists at each (proper) time (of the moving particle) a (local) rest frame where the acceleration of the particle is Brownian in any spacelike direction of the frame, when computed using the time of the rest frame. The processes of this class were named \textit{relativistic diffusions} in reference to the diffusion phenomenon they modelize. The authors of the article have started the study of this class developing an analytical approach to the situation based on a transport equation. We would like to propose in the present article a pathwise approach to this class of processes on a general Lorentzian manifold. With in mind the diffusion phenomenon of colloidal particles in fluids, we shall describe their dynamics as random perturbations of differential equations. In the spirit of the work of Franchi and Le Jan, we shall lift these dynamics to the frame bundle of the manifold, where they will be defined as flows of stochastic differential equations. This framework will enable us to re-prove directly many of the results obtained so far as well as new results and prospects stemming from the pathwise nature of our approach. 

\smallskip

We have organized the exposition as follows. Section \ref{Section1} is dedicated to describing the class of relativistic processes in Minkowski spacetime, so as to separate probability and geometry problems. The class of relativistic diffusions is thus motivated and defined in section \ref{ModelsInMin}. We give in section \ref{SectionOnePart} a probabilistic definition of the one-particle distribution function for each relativistic diffusion, and prove that it satisfies a fundamental equation. Section \ref{DiffLorMan} is dedicated to investigating the general situation where the geometric background is any Lorentzian manifold. After having defined the dynamics in the orthonormal frame bundle in section \ref{SubsubsectionDefnDiff}, we shall spend some time  in section \ref{SubsubsectionSubDiff} looking at what can happen in the unit sub-bundle of the tangent bundle. We shall define in section \ref{SubSubSectionOnePartFunction} the one-particle distribution function for each relativistic diffusion under a mild hypothesis on the global geometry of spacetime. The relevance of this notion in the study of the Poisson and Martin boundaries of the Franchi-Le Jan process will be discussed in section \ref{SubSubSectionHarmonicFunctions}. Finally, we shall prove in section \ref{BouquetFinal} a general $H$-theorem. A number of open problems are scattered throughout the text. Numerous examples have been included so as to help the reader to get an idea of the state of the field.

\medskip

\noindent \textbf{Notation.} We shall write ${\circ d}$ for the Stratonovich differential. The sign $d$ will be used for the usual differentiation with respect to the time, or for Ito's differential.


\section{Relativistic diffusions in Minkowski spacetime}
\label{Section1}

\subsection{Definitions and examples}
\label{ModelsInMin}

\paragraph{a) Geometric framework.} Recall Minkowski space is the product $\RR\times\RR^3$ equipped with the metric
$$
\forall\,\zeta=(t,x)\in\RR^1\times\RR^3,\quad q(\zeta) = t^2-\left(\bigl(x^1\bigr)^2+\bigl(x^2\bigr)^2+\bigl(x^3\bigr)^2\right),
$$
if we write $(t,x^1,x^2,x^3)$ for the co-ordinates of $\zeta$ in the canonical basis $\bigl\{\varep^0,\varep^1,\varep^2,\varep^3\bigr\}$ of $\RR\times\RR^3$. To distinguish Minkowski spacetime from the Euclidean space $\RR^4$, we shall denote the former by $\RR^{1,3}$. The half-unit sphere
$$
\HH = \{\zeta=(t,x)\in\RR^{1,3}\,;\,q(\zeta)=1,\,t>0\}
$$
inherits from the ambient (non-definite positive) metric $q$ a Riemannian metric of constant curvature, which makes it a model of the ($3$-dimensional) hyperbolic space. As any $\mcC^1$ timelike path can be re-parametrized so that its speed should belong to $\HH$, we shall look at the space $\RR^{1,3}\times\HH$ as the configuration space of timelike $\mcC^1$ trajectories of a point of $\RR^{1,3}$. The set of direct linear isometries of $q$ is the group $SO(1,3)$. Any element $\bfg$ of $SO(1,3)$ represents a rest frame $\bfg = \bigl(\bfg^0,\bfg^1,\bfg^2,\bfg^3\bigr)$ of $\RR^{1,3}$. The function $\zeta\in\RR^{1,3} \mapsto q(\bfg^0,\zeta)$ will be called the \textit{\textbf{time function associated with the frame $\bfg$}}.

\smallskip

It will also be fruitful to define the motion of a(n infinitesimally small) rigid object. The configuration space of this dynamics will be the set $\RR^{1,3}\times SO(1,3)$. We shall look at a point $\bigl(m,\bigl(\bfg^0,\bfg^1,\bfg^2,\bfg^3\bigr)\bigr)$ as the infinitesimal rigid object\footnote{$\delta$ is some infinitesimal positive number.} $m+\textrm{ConvHull}(\delta\bfg^1,\delta\bfg^2,\delta\bfg^3)$ contained in the affine spacelike hyperplane $m+span\bigl(\bfg^1,\bfg^2,\bfg^3\bigr)$, and having $4$-velocity $\bfg^0$. An element of $\RR^{1,3}\times SO(1,3)$ can also be seen as an observer.

Notice that $SO(1,3)$ has $4$ connected components; we shall denote by $SO_0(1,3)$ the connected component of the identity. To shorten notations, we shall write $\OO\RR^{1,3}$ for $\RR^{1,3}\times SO_0(1,3)$.

The introduction of the following notations will clarify the description of the dynamics we are interested in. We shall denote by $E_i\in so(1,3)$ the Lie element such that $\exp(tE_i)$ is the hyperbolic rotation of angle $t$ in the $2$-dimensional plane generated by $\varep^0$ and $\varep^i$. In matrix notations
$$
E_1=\begin{pmatrix} 0 & 1 & 0 & 0 \\ 1 & 0 & \cdots \\ 0 & \vdots & & \\ 0 & \vdots & {\bf O}_2\end{pmatrix},\; E_2=\begin{pmatrix} 0 & 0 & 1 & 0 \\ 0 & 0 &  & \cdots \\ 1 & \vdots &  &  \\ 0 & \vdots & & {\bf O}_2\end{pmatrix},\; E_3=\begin{pmatrix} 0 & 0 & 0 & 1 \\ 0 & 0 &  & \cdots \\ 0 & \vdots &  &  \\ 1 & \vdots & & {\bf O}_2\end{pmatrix}. 
$$
Four vector fields on $\OO\RR^{1,3}$ will be of particular interest.
\begin{equation}
\label{CanonicalHorVectFields}
\begin{split}
& H_0\bigl((m,\bfg)\bigr) = (\bfg^0,0), \\
& \textrm{for } i=1..3, \quad V_i\bigl((m,\bfg)\bigr) = (0,gE_i).
\end{split}
\end{equation}

Note that the $\RR^{1,3}$-part of the integral lines of the vector field $H_0$ are the geodesics of $\RR^{1,3}$, which are straight lines. We shall set $\HH_m\RR^{1,3} = \bigl\{(m,\zeta)\in\RR^{1,3}\times\HH\bigr\}$ and write $\OO_m\RR^{1,3}$ for $\bigl\{(m,\bfg)\in\OO\RR^{1,3}\,;\,\bfg\in SO_0(1,3)\bigr\}$. 

\medskip

An important feature of our approach to relativistic diffusions is that we have chosen to describe the dynamics in the phase space $\OO\RR^{1,3}$, where it has a natural and simple form; this corresponds to look at the motion of a small rigid object. We shall look at what happens in $\RR^{1,3}\times\HH$ in a later section.

\paragraph{b) Dynamics.}
$\bullet$ \textbf{Unperturbed system.}   We have indicated in the introduction that relativistic diffusions should be considered as a class of toy models of diffusion in different media. We are going to define them as random perturbations of deterministic evolutions given by the flow of a vector fields $V$ on $\RR^{1,3}$. With in mind diffusion of particles in a fluid, we shall make the hypothesis that  $V$ has no $\RR^{1,3}$-part and acts only on the $SO(1,3)$-part of $\OO\RR^{1,3}$, although this assumption could be relaxed. The unperturbed sytsem is defined by the differential equation 
\begin{equation}
\label{UnperturbedSystem}
\begin{split}
& dm_s = \bfg^0_s\,ds, \\
& d\bfg_s = V(\bfg_s)ds.
\end{split}
\end{equation}
Note that the requirement that $\frac{dm_s}{ds} = \bfg^0_s\in\HH$ implies that the parameter $s$ is the proper time of the timelike path $\{m_s\}_{s\geq 0}$ of $\RR^{1,3}$.

\smallskip

$\bullet$ \textbf{Action of the surrounding medium.} How should we model the form taken by the random perturbation of the dynamics associated with a given medium? Maybe the proper way to proceed would consist in giving first a description of the microscopic thermodynamical and electro-magnetical properties of the medium in order to put forwards the source of randomness, and to infer from this description a description of the random perturbation it induces on the dynamics of a test object. We have chosen to propose a rather general action model which should convey the essential features of many situations, and not to model the medium itself.

\smallskip

The action of the fluid on the moving object $\{\be_s\}_{s\geq 0} = \bigl\{(m_s,\bfg_s)\bigr\}_{s\geq 0}$ will be represented by the datum of an $\OO\RR^{1,3}$-valued \textit{previsible} process $\{\ff_s\}_{s\geq 0}$ such that $\ff_s(\be_.) = \ff_s\bigl((m_.,\bfg_.)\bigr) = (m_s,f_s)$ for some orthonormal basis $f_s = \left(f^0_s,f^1_s,f^2_s,f^3_s\right)$ of $T_{m_s}\RR^{1,3}$(\footnote{Note that $\ff_s$ and $\be_s$ have the same $\RR^{1,3}$-part equal to $m_s$.}). The random perturbation induced by the medium on the dynamics results in adding to the deterministic acceleration a random part which is determined by the following requirement. \textit{When computed in the rest frame $\ff_s$, i.e. using its associated time}, the acceleration of $m_s$ has a deterministic part and a random part which is Brownian in any spacelike direction belonging to $\textrm{span}(f^1_s,f^2_s,f^3_s)$. To complete this description, we shall ask the vectors $\bfg^1_s,\bfg^2_s,\bfg^3_s$ to be transported parallelly along the "Brownian" increment of $\bfg^0_s$.

\paragraph{c) A preliminary example.}
Before giving a mathematically clean definition of this class of processes, let us look at the heuristic description of what happens when $V=0$ and the 'vertical' action process $\ff_.$ is constant, equal to Id, \textit{i.e.} $f_s=\{\varep^0,...,\varep^3\}$ for any $s$.

\smallskip

Denote by $\bigl\{(m_s,\bfg_s)\bigr\}_{s\geq 0}$ the $\OO\RR^{1,3}$-valued process corresponding to these data and write $t_s$ for the $\varep^0$-component of $m_s$. As we have $dm_s = \bfg^0_s\,ds$, the function $s\mapsto t_s$ is a $\mcC^1$ increasing function that can be used as a parameter of the process. Given $t\in\RR$, set $\tau_t = \inf \bigl\{s\geq 0\,;\,t_s=t\bigr\}$ and look at the re-parametrized process $\bigl\{(m_{\tau_t},\bfg_{\tau_t})\bigr\}_{t\geq q(\varep^0,m_0)}$; denote it by $\bigl\{(\widehat m_t,\widehat\bfg_t)\bigr\}_{t\geq q(\varep^0,m_0)}$. The above description of the action of the surrounding medium on the dynamics means that the $\textrm{span}(\varep^1,\varep^2,\varep^3)$-part of $\displaystyle{d\widehat\bfg^0_t}$ is a Brownian increment.

\begin{figure}[h!]
\label{FigROUP}
\begin{center}
\input{FigureROUP.pstex_t} 
\end{center}
\caption{Dynamics when $\ff=\textrm{Id}$ and $V=0$}
\end{figure}

The Brownian spacelike part $\displaystyle{\sum_{i=1..3}}\varep^i{\circ d\widehat w^i_t}$ of the increment of the speed can be seen in figure 1, in red; the increment itself is in green. The notation $\widehat w$ stands here for a $3$-dimensional Brownian motion. If we write $\circ d\widehat\bfg^0_t = \displaystyle{\sum_{i=1..3}}\widehat\bfg^j_t{\circ d\widehat\beta^j_t}$, then
$$
\circ d\widehat w^i_t = -\sum_{j=1..3} q(\varep^i,\widehat\bfg^j_t)\,{\circ d\widehat\beta^j_t}.
$$
Denote by $A(\bfg)$ the $3\times 3$ matrix with coefficients $(i,j)\in [1,3]^2$ equal to $q(\varep^i,\bfg^j)$. This matrix being invertible,
\begin{equation}
\label{IdentityBeta}
{\circ d\widehat\beta_t} = -A(\widehat\bfg_t)^{-1}{\circ d\widehat w_t}.
\end{equation}
Back to the proper time $s$ of the process, we shall write $\displaystyle{{\circ d}\bfg^0_s = \sum_{j=1..3}\bfg^j_s\,{\circ d}\beta_s^j}$. Write $A_s$ for $A(\bfg_s)$. Identity \eqref{IdentityBeta} implies that
$$
{\circ d\beta_s} = q(\varep^0,\bfg^0_s)^{\frac{1}{2}}\,A_s^{-1}{\circ dw_s}
$$
for some $3$-dimensional \bm $w$. The $\RR^3$-valued process $\beta$ is the process that really drives the dynamics. Last, we shall ask the vectors $\bfg^1_s,\bfg^2_s,\bfg^3_s$ to be parallelly transported along the paths $\{\bfg^0_s\}_{s\geq 0}$ in $\HH$. The above heuristic description gives rise to the following equations of motion
\begin{equation*}
\begin{split}
&{\circ d} m_s =  \bfg^0_s\,ds, \\
&{\circ d}\bfg_s = \bfg_s E_i\,{\circ d\beta^i_s}.
\end{split}
\end{equation*}

\paragraph{d) Definition.}
We shall now leave appart this example to write down the equations of the dynamics of $\bigl\{(m_s,\bfg_s)\bigr\}_{s\geq 0}$ corresponding to general data $V$ and $\ff$. Recall the surrounding medium will be represented by the datum of a previsible process $\{\ff_s\}_{s\geq 0} = \{\ff_s(\be_.)\}_{s\geq 0}$ \st $\ff_s = (m_s,f_s) = \bigl(m_s,(f^0_s,...,f^3_s)\bigr)$ belongs to $\OO_{m_s}\RR^{1,3}$. Its action on the dynamics has been heuristically described in paragraph \textbf{b)}. Define the random matrix process $\{A_s\}_{s\geq 0}$, with coefficient $(i,j)\in [1,3]^2$ equal to $q(f_s^i,\bfg_s^j)$ at time $s$; set
\begin{equation}
\label{DefnBeta0}
{\circ d\beta_s} = q(f_s^0,\bfg^0_s)^{\frac{1}{2}}\,A_s^{-1}{\circ dw_s}.
\end{equation}

\begin{defn}
\label{DefnRelDiffMin}
Define the $\RR^3$-valued process $\beta$ as above. A \textbf{$(V,\ff)$-diffusion} is a process $\{\be_s\}_{s\geq 0} = \bigl\{(m_s,\bfg_s)\bigr\}_{s\geq 0}$ satysfying the stochastic differential equations
\begin{equation}
\label{EqDyn1}
\begin{split}
&{\circ d} m_s =  \bfg^0_s\,ds \\
&{\circ d}\bfg_s = V(\be_s)ds + \bfg_sE_i\,{\circ d\beta^i_s},
\end{split}
\end{equation}
where Einstein's summation convention is used, as in the sequel.
\end{defn}
Using notations \eqref{CanonicalHorVectFields}, equation \eqref{EqDyn1} can be written
\begin{equation}
\label{SDEDynamics}
{\circ d}\be_s = H_0(\be_s)ds + V(\be_s)ds + V_i(\be_s)\,{\circ d\beta_s^i}.
\end{equation}
In reference to the interpretation of $\OO\RR^{1,3}$ in terms of infinitesimal rigid objects given in paragraph \textbf{a)}, this equation can be interpreted as describing the random motion of an infinitesimal rigid object in $\RR^{1,3}$; there are nonetheless no need to understand it that way if you do not feel comfortable with infinitesimals. In any case, the simple and intrinsic character of this equation should be compared with the co-ordinate approach proposed up to now, as presented for instance in the article \cite{ROUPUnifying1} of C. Chevalier and F. Debbasch. The simplicity of the formalism of stochastic differential equations will enable us not to rely on the covariant treatment used so far.

Note that since $\ff_s(\be_.)$ might depend on the whole history of $\be_.$ until time $s$, the increment ${\circ d\beta_s}$ shares this property, and equations \eqref{DefnBeta0} and \eqref{SDEDynamics} do not generally define a Markov process. This might be relevant from a modelization point of view if we consider an object having internal parameters evolving with time, and whose value at proper time $s$ could influence the way the surrounding medium acts on it. Let us give three (Markovian) examples before commenting any further.

\paragraph{e) Previously studied examples.} 
Three $(V,\ff)$-diffusions have attracted attention up to now. 
\begin{enumerate}
    \item The Dudley(-Franchi-Le Jan) process introduced by Dudley in \cite{Dudley1} (and generalized in \cite{FranchiLeJan} by Franchi and Le Jan) is a perturbation of the geodesic flow. It corresponds to taking $V=0$ and $\ff_s=\be_s$. The dynamics driving process $\beta$ is then equal to the \bm $w$, and no time-change is needed\footnote{That is, the time scaling $q(f^0_s,\bfg^0_s)$ is here equal to $1$.}. It is described in a simple way saying that 
\begin{itemize}
   \item $\{\bfg^0_s\}_{s\geq 0}$ is a \bm on the hyperbolic space $\HH$,
   \item $(\bfg^1_s,\bfg^2_s,\bfg^3_s)\in T_{\bfg^0_s}\HH$ is obtained from $(\bfg^1_0,\bfg^2_0,\bfg^3_0)$ by parallel transport along the path $\{\bfg^0_r\}_{0\leq r\leq s}$, and
   \item $m_s = m_0 + \int_0^s \bfg_r^0\,dr$.
\end{itemize}
This process is the only process determined entirely by the datum of the geometric background (a result due to Dudley in \cite{Dudley1}(\footnote{Note that we have uniqueness up to a time scaling by a constant in the $\HH$-\bm $\{\bfg^0_s\}_{s\geq 0}$.})). This property gives it a special position in the family of $(V,\ff)$-diffusions. Yet, its drawback as a model in Minkowski spacetime of a diffusing particle is that, except if we locate the source of motion in the particle itself, it is not clear what entity could give rise to such an interaction process. So it might be less interesting from a modelization point of view. Consult yet the article \cite{DowkerHensonSorkin} of Dowker, Henson and Sorkin for a physical motivation from quantum mechanics. Nevertheless, the long-time behaviour of this process and its Lorentzian version may have many things to say about the geometry at infinity of spacetime; this might happen to be of some (theoretical) physical interest. We shall discuss this point in section \ref{SubSubSectionHarmonicFunctions}.
   
   \item The relativistic Ornstein-Uhlenbeck process (R.O.U.P.) was introduced by F. Debbasch, K. Mallick and J.P. Rivet in the article \cite{ROUP0} as a model of diffusing particle in a fluid at equilibrium. It corresponds to the $\bigl(V,\textrm{Id}\bigr)$-diffusion with
$$
V\bigl((m,\bfg)\bigr) = -\al\,\textrm{grad}(\ln \gamma)
$$
for some positive constant $\al$.  We have written here $\gamma$ for $q(\varep^0,\bfg^0)$ and grad for the gradient in $\HH$. In this case, the dynamics driving process $\beta$ is not equal to the \bm $w$. The existence for this process of an invariant measure of the form\footnote{The measure $d\bfg$ is a Haar measure on the unimodular group $SO_0(1,3)$, and $a$ and $b$ are positive constants.} $ae^{-b\gamma}dm\otimes d\bfg$ found by J\"uttner in \cite{Juttner} was a motivation for its introduction; see the introduction of the article \cite{ROUP0}. We shall see in the general framework of section \ref{SubsubsectionSubDiff}, that this $\OO\RR^{1,3}$-valued diffusion gives rise to an $\HH\RR^{1,3}$-valued diffusion, which is the R.O.U.P. as defined in \cite{ROUP0} and the subsequent works of the authors and their co-authors.

   \item Last, Dunkel and H\"anggi introduced in their article \cite{DunkelHanggi} a kind of mixing of the previous two models in which the frame $\ff_s=\be_s$, as in the Dudley-Franchi-Le Jan diffusion, and $V$ is constructed in such a way that the process admits the same awaited invariant measure as the R.O.U.P. 
\end{enumerate}
We shall come back to these models in the general framework of section \ref{DiffLorMan}.

\paragraph{f) Non-isotropic medium.} 
This way of defining $(V,\ff)$-diffusions has the advantage to be flexible enough to provide models of what should be a relativistic diffusion in a non-isotropic medium. We shall take into account the non-isotropy of the motion replacing the up to now isotropic input ${\circ dw_s}$ by a non-isotropic semimartingale in equations \eqref{DefnBeta0} and \eqref{EqDyn1} of dynamics. Setting for instance $M = \textrm{diag}(1,1,2)$ and denoting by $\{B_s\}_{s\geq 0}$ an $\RR^3$-valued Brownian motion, the use in the R.O.U.P. dynamics of an input ${\circ dw_s} = M{\circ dB_s}$ will give rise to a motion in a medium where one spacelike (fixed) direction differs from the others. One could also replace $w$ by any continuous semi-martingale to adapt the model to a given situation. Jumps could even be introduced to take into account possible shocks.

The article \cite{JurgenFranchi} of J. Franchi and J. Angst proposes another model in Minkowski spacetime of random dynamics in a non-isotropic medium.

\paragraph{g)  Probabilistic matters.}
Let us be more precise in the definition of a $(V,\ff)$-diffusion\footnote{Refer to the chapter $V.8$ of the book \cite{RogersWilliams} by Rogers and Williams for all this paragraph.}. Let $\bigl(W,\{\mcH_t\}_{t\geq 0}\bigr)$ denote the Polish space $\mcC\bigl(\RR_+,\OO\RR^{1,3}\bigr)$, endowed with the filtration generated by its co-ordinate process. Let $\ff : \RR_+\times W \rightarrow \OO\RR^{1,3}$ be a previsible path functional. A $(V,\ff)$-diffusion will consist in the datum of a filtered probability space $\bigl(\Omega,\{\mcF_t\}_{t\geq 0},\PP\bigr)$ satisfying the usual conditions, an $\bigl(\{\mcF_t\}_{t\geq 0},\PP\bigr)$-\bm $w$ on $\RR^3$, and a $\mcC\bigl(\RR_+,\OO\RR^{1,3}\bigr)$-valued process $\be$ defined on $\bigl(\Omega,\{\mcF_t\}_{t\geq 0}\bigr)$ such that equations \eqref{DefnBeta0} and \eqref{EqDyn1} hold. These sorts of details will be implicit in the sequel.

Existence and uniqueness results exist for equations such as \eqref{DefnBeta0} and \eqref{EqDyn1}. Consult \cite{RogersWilliams} and the references given therein for example. These issues will raise no problem in the example we shall consider. 

We should apologize for the mis-use of the word "diffusion" in this context, as it is usually used when $\ff_s(\be_.) = \ff(\be_s)$, which is not supposed here. We have chosen to keep this denomination in reference to the situation it modelizes. The word "diffusion" will keep in the sequel its usual meaning, and we shall always write $(V,\ff)$-diffusion for a process of our class. 

Last, we shall use the notation $\{\be_s\}_{s\geq 0}$, indexing the trajectories by $\RR_+$, regardless of the possibly finite lifetime of the process. One can add a cemetery point to the space to deal with such issues.

\subsection{One-particle distribution function of Markovian $(V,\ff)$-diffusions}
\label{SectionOnePart}

As explained in the introduction, the main aim of his article is to convince the reader of the usefulness of a pathwise approach to relativistic diffusions. This section will illustrate this point giving a clear definition of the one-particle distribution function of a $(V,\ff)$-diffusion. We refer to the article \cite{DRVanLeeuwen} of F. Debbasch, J.P. Rivet and W.A. van Leeuwen for a physical discussion of this concept of statistical physics and for the interest of a clear definition of this notion\footnote{The article \cite{IsraelBoltzmann} of W. Israel can also be consulted on this subject.}. We shall investigate the general situation on a Lorentzian manifold in section \ref{PreBouquetFinal}. Let us first describe the framework of the problem.

\paragraph{a) Framework.} 
We shall suppose in this section that 
$$
\ff_s(\be_.) = \ff(\be_s)
$$
for some function $\ff : \OO\RR^{1,3}\rightarrow\OO\RR^{1,3}$ \st $\ff\bigl((m,\bfg)\bigr) = \Bigl(m,\bigl(f^0(\be),...,f^3(\be)\bigr)\Bigr)$. It follows that the process $\{\be_s\}_{s\geq 0} = \bigl\{(m_s,\bfg_s)\bigr\}_{s\geq 0}$ is an $\OO\RR^{1,3}$-valued Markov process. Write $A(\be)$, or simply $A$, for the $3\times 3$ matrix with coefficient $(i,j)\in [1,3]^2$ equal to $q(f^i(\be),\bfg^j)$. The generator of the $(V,\ff)$-process is given by the formula
\begin{equation}
\label{Generator}
L:= H_0 + V + \frac{\la}{2}\,V_iB^{ij}V_j,
\end{equation}
where $B = \left(A^{-1}\right)^*A^{-1}$ is a $3\times 3$ non-negative symmetric matrix. Here as in the sequel, a vector field is seen as a first order differential operator; so, an expression like $V_iB^{ij}V_j f$ should be more properly written $V_i\bigl(B^{ij}V_j(f)\bigr)$. The use of the notation 
$$
\la := q(f^0(\be),\bfg^0)
$$ 
will be useful to shorten formulas, here as in the sequel. Recall we have supposed that the flow of $V$ preserves each fiber of the projection $(m,\bfg)\in\OO\RR^{1,3}\rightarrow m\in\RR^{1,3}$.

We shall denote by $d\bfg$ the Haar measure on (the unimodular group) $SO_0(1,3)$ whose image by the projection $\bfg\in SO_0(1,3)\mapsto \bfg^0$ is the Riemannian measure on $\HH$. Last, we shall associate to any subset $A$ of $\RR^{1,3}$ the (principal) bundle
$$
\OO A :=\bigl\{(m',\bfg')\in\OO\RR^{1,3}\,;\,m'\in A,\,\bfg'\in SO_0(1,3)\bigr\}.
$$
If $A$ is a spacelike hypersurface of $\RR^{1,3}$, denote by $\sigma_A(dm')$ the volume measure induced by $q$ on $A$; we define the measure
$$
\Vol_{\OO A}(d\bfg'\wedge dm') := d\bfg'\otimes\sigma_A(dm')
$$
on the bundle $\OO A$.

\paragraph{b) One-particle distribution function.} 
A few more notations are needed to define the one-particle distribution function and state its main properties. Fix a point $\be=(m,\bfg)\in\OO\RR^{1,3}$, and define the hyperplane of $\RR^{1,3}$
$$
\VV_\be = \bigl\{m'\in\RR^{1,3}\,;\,m'\in m+\bigl(\bfg^0\bigr)^\perp\bigr\}; 
$$
denote by $H_\be$ the hitting time 
$$
H_\be = \inf\{s>0\,;\,\be_s\in \OO\VV_\be\}.
$$
We shall associate to any $\al\in SO_0(1,3)$ and any $t\in\RR$ the hyperplane $\VV_{t,\al} := \{m'\in\RR^{1,3}\,;\,q(m',\al^0)=t\}$ and the hitting time $H_{t,\al} = \inf\{s>0\,;\,\be_s\in\OO\VV_{t,\al}\}$.

Notice that the Liouville measure induced by $q$ on $\OO\RR^{1,3}$ is the product measure 
$$
\Vol(d\bfg\wedge dm) := d\bfg \otimes \LEB_4(dm).
$$ 
We shall denote by $L^*$ the $\LL^2(\Vol)$-dual of the operator $L$; we have $V_i^* = -V_i$ and $H_0^* = -H_0$.

\begin{ThmDefn}
\label{PropDefnOnePartDistrFunction}
\begin{enumerate}
   \item Let $\be_0\in\OO\RR^{1,3}$ be different from $\be$. The \rv $\be_{H_\be}{\bf 1}_{H_\be<\infty}$ has under $\PP_{\be_0}$ a smooth density $f\bigl(\be_0\,;\,(m',\bfg')\bigr)$ \wrt the measure $\Vol_{\OO\VV_\be}(d\bfg'\wedge dm')$ on $\OO\VV_\be$.
\end{enumerate}
\noindent The function $\be\in\OO\RR^{1,3}\backslash\{\be_0\}\mapsto f(\be_0\,;\,\be)$ is called the \emph{\textbf{one-particle distribution function of the $(V,\ff)$-diffusion started from $\be_0$}}.
\begin{enumerate}
   \item[\textit{2.}] We have
\begin{equation}
\label{InterpretationOnePartMin}
\EE_{\be_0}\big[f(\be_{H_{t,\al}})\bigr] = \int f(\be)\,q\bigl(\al^0,\bfg^0\bigr)f(\be_0\,;\,\be)\,\Vol_{\OO\VV_{t,\al}}(d\be)
\end{equation}
for any bounded function $f$ on $\OO\VV_{t,\al}$.
   
   \item[\textit{3.}] The function $f(\be_0\,;\,\cdot)$ satisfies the equation
\begin{equation}
\label{FundamentalEq}
L^*f(\be_0\,;\,\cdot) = 0
\end{equation}
on $\OO\RR^{1,3}\backslash\{\be_0\}$.
\end{enumerate}
\end{ThmDefn}

It is clear from its definition that this function is defined in an intrinsic way; physicists use to say that $f(\be_0\,;\,\be)$ is a \textit{Lorentz scalar}. 
We shall prove in section \ref{SubSubSectionOnePartFunction} a similar theorem in the general framework presented in section \ref{DiffLorMan}. We have chosen to present here a heuristic proof of point \textit{2} and to give a detailed proof of the general statement after proposition \ref{InterpretationOnePart}, in section \ref{SubSubSectionOnePartFunction}. Points \textit{1} and \textit{3} of theorem/definition \ref{PropDefnOnePartDistrFunction} are proved in detail below.

\begin{Dem}
\textit{1.} The strategy of the proof is simple. Given $\be=(m,\bfg)$, we are going to re-parametrize the process as a function of the time associated to the frame $\bfg$ and see that $f(\be_0\,;\,\cdot)$ is the density \wrt $\Vol_{\OO\VV_\be}$ of the position at some fixed time of a hypoelliptic diffusion.

\begin{figure}[h!]
\label{FigROUP}
\begin{center}
\input{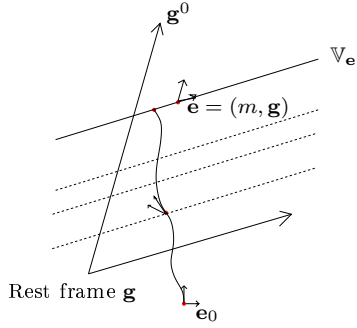} 
\end{center}
\caption{Re-parametrized process}
\end{figure}

Define the \textit{chronological past of $\OO\VV_\be$} as the set $I^-(\OO\VV_\be)$:
\begin{equation*}
\begin{split}
\Bigl\{\bigl(\gamma(0),\bfg'\bigr)\in\OO\RR^{1,3}\,;\,\gamma \textrm{ future-oriented timelike path from }\gamma(0) \textrm{ to a point of} &\textrm{ the set } m+\bigl(\bfg^0\bigr)^\perp, \\
&\bfg'\in SO(1,3)\Bigr\}.
\end{split}
\end{equation*}
The \rv $\be_{H_\be}{\bf 1}_{H_\be<\infty}$ being identically equal to $0$ if $\be_0$ does not belong to the chronological past of $\OO\VV_\be$, we shall suppose in the sequel that $\be_0$ belongs to it, in which case $H_\be$ is almost surely finite.

\smallskip

Set $t_0^{\bfg} = q(m_0,\bfg)$ and define the stopping times 
$$
\forall\,t\in \RR, \quad S^\bfg_t = \inf\{s\geq 0\,;\,q(\bfg^0,m_s) = t\}.
$$
The process $\{\be_{S_t}\}_{t \geq t_0^{\bfg}}$ is the process $\be_.$ re-parametrized by the time associated with $\bfg$. It has generator 
$$
\frac{1}{\la_\bfg}L,
$$
where $\la_\bfg = \la_\bfg(\bfg') = q(\bfg^0,\bfg'^0)$. We shall write
$$
\be_{S^\bfg_t} = \bigl((t\bfg^0+x_{S_t^\bfg}),\bfg'_{S^\bfg_r}\bigr)\in \OO\RR^{1,3}, \quad\textrm{with }x_{S_t^\bfg}\in \textrm{span}(\bfg^1,\bfg^2,\bfg^3),
$$
and shall look at 
$$
\be^{\bfg}_t := (x_{S_t^\bfg},\bfg'_{S^\bfg_t}) \in \OO\,\textrm{span}(\bfg^1,\bfg^2,\bfg^3).
$$
The \rv $\be_{H_\be}$ is equal to $m_0+t_1\bfg^0 + \be_{t_1}^\bfg$, with $t_1=q(m-m_0,\bfg^0)$. We shall prove the first point of theorem \ref{PropDefnOnePartDistrFunction} showing that the $\OO\textrm{span}\bigl(\bfg^1,\bfg^2,\bfg^3\bigr)$-valued diffusion $\be^\bfg_.$ is a hypoelliptic diffusion. The distribution at time $t_1$ of this diffusion will then have a smooth density \wrt the volume element on $\OO\textrm{span}\bigl(\bfg^1,\bfg^2,\bfg^3\bigr)$ to be defined below. The measure $\Vol_{\OO\VV_\be}$ being the image of the volume element by $\RR^{1,3}$-translation by $m_0+t_1\bfg^0$, this will imply that $\be_{H_\be}$ has a smooth density \wrt $\Vol_{\OO\VV_\be}$.

\smallskip

To complete this program we shall denote by $\ba^\bfg:=(x,\bfg')$ a generic element of $\OO\,\textrm{span}(\bfg^1,\bfg^2,\bfg^3)$. Note that since $\bfg$ has determinant equal to $1$, the change of variable formula says us that the volume element induced by $q$ on the $3$-dimensional vector space spanned by $\bfg^1,\bfg^2$, and $\bfg^3$ is the $3$-dimensional Lebesgue measure. We shall write $\Vol_{\OO\,\textrm{sp}(\bfg^1,\bfg^2,\bfg^3)}(d\bfg'\wedge dx) = d\bfg'\otimes\LEB_3(dx)$ the volume measure on the bundle $\OO\,\textrm{span}(\bfg^1,\bfg^2,\bfg^3)$.

\smallskip

To describe the generator $L^\bfg$ of the process $\be^{\bfg}_.$, denote first by $\partial_x$ the differentiation operation in the direction of the vector space $\textrm{span}(\bfg^1,\bfg^2,\bfg^3)$ and decompose $\bfg'^0$ as
$$
\bfg'^0 = \la_\bfg\bfg^0 + \sum_{i=1}^3\dot{x}_i\bfg^i.
$$
Using these notations, we can write for any smooth function $f$
$$
\frac{Lf}{\la_\bfg} = \partial_t f + \frac{\bigl(\partial_x f\bigr)(\dot{x})}{\la_\bfg} + \frac{Vf}{\la_\bfg} + \frac{\la}{2\la_\bfg}\,V_iB^{ij}V_j f.
$$
So the generator $L^\bfg$ of the process $\be^\bfg_.$ is given by the formula
$$
L^\bfg f = \frac{\partial_x f(\dot{x})}{\la_\bfg} + \frac{Vf}{\la_\bfg} + \frac{\la}{2\la_\bfg}\,V_iB^{ij}V_j f.
$$
Write $h_{\be_0}^\bfg(t,d\ba^\bfg)$ for the law of $\be^{\bfg}_t, t\geq t_0^\bfg$. As is well known, these distributions satisfy the heat equation 
\begin{equation}
\label{HeatEqAlpha}
\partial_t\,h^\bfg_{\be_0} = \bigl(L^\bfg\bigr)^{*\bfg} h^\bfg_{\be_0},
\end{equation}
where $\bigl(L^\bfg\bigr)^{*\bfg}$ is the $\LL^2\bigl(\Vol_{\OO\,\textrm{sp}(\bfg^1,\bfg^2,\bfg^3)}\bigr)$-dual of the operator $L^\bfg$. Since the matrix $B = \bigl(A^{-1}\bigr)^*A^{-1}$ is symmetric and $V_i^{*\bfg} = -V_i$, we have  
$$
\bigl(L^\bfg\bigr)^{*\bfg} h^\bfg_{\be_0}= -\left(\partial_x h^\bfg_{\be_0}\right)\Bigl(\frac{\dot{x}}{\la_\bfg}\Bigr) + V^{*\bfg} \Bigl(\frac{h^\bfg_{\be_0}}{\la_\bfg}\Bigr) + \frac{1}{2}\bigl(V_iB^{ij}V_j\bigr)\Bigl(\la\frac{h^\bfg_{\be_0}}{\la_\bfg}\Bigr).
$$
It is easy to see on this formula that the operator $\partial_t - \bigl(L^\bfg\bigr)^{*\bfg}$ on $\RR\times\OO\,\textrm{span}(\bfg^1,\bfg^2,\bfg^3)$ satisfies H\"ormander's criterium for hypoellipticity. It follows that $h_{\be_0}^\bfg(t,\cdot)$ has a smooth density \wrt the measure $\Vol_{\OO\,\textrm{sp}(\bfg^1,\bfg^2,\bfg^3)}(d\ba^\bfg)$ on $\OO\,\textrm{span}(\bfg^1,\bfg^2,\bfg^3)$, for $t>q(m_0,\bfg^0)$. We have seen that it implies that $\be_{H_\be}$ has a smooth density \wrt $\Vol_{\OO\VV_\be}$.

\bigskip

\noindent \textit{2.} As said above, we present here a heuristic proof of point \textit{2}. The reader will find the detailed proof of the general statement after proposition \ref{InterpretationOnePart}, in section \ref{SubSubSectionOnePartFunction}. We are going to explain the situation for Dudley's process, nothing else than additionnal notations being necessary to understand the general case of Markovian $(V,\ff)$-diffusions.

\smallskip

$\bullet$ We shall get a clearer image of the situation considering the continuous dynamics described by equation 
\begin{equation}
\label{EqDudleyDiff}
d\be_s = H_0ds + V_i{\circ dw^i_s}
\end{equation}
as the dynamics of a \rw $\{\widetilde\be_s\}_{s\geq 0} = \bigl\{(\widetilde m_s,\widetilde\bfg_s)\bigr\}_{s\geq 0}$ making infinitesimal steps.

Given an 'infinite' integer $N$ (\textit{i.e.} a nonstandard hyperfinite integer), the quantity $\frac{1}{N}$ is a positive infinitesimal. Let us denote by $\{\Delta_k\}$ a 'sequence' of iid $\RR^d$-valued centered Gaussian \rvs with variance $\frac{1}{N}$. The dynamics of the \rw is defined on each interval of the form $\bigl[\frac{k}{N},\frac{k+1}{N}\bigr), k\geq 1$ as follows.
\begin{itemize}
   \item The process $\widetilde\bfg_s$ has a jump at time $\frac{k}{N}$: $\displaystyle{\quad \widetilde\bfg_{\frac{k}{N}} = \widetilde\bfg_{\left(\frac{k}{N}\right)^-}.\exp\bigl(E_i\,\Delta w^i_k\bigr)}$. The process $\{\widetilde m_s\}_{s\geq 0}$ has no jumps at that time.
   \item $\widetilde\bfg_s$ is constant and $d\widetilde m_s = \widetilde\bfg^0_sds$, in the time interval $\bigl(\frac{k}{N},\frac{k+1}{N}\bigr)$.
\end{itemize}
$\widetilde\bfg_s$ is constant and $d\widetilde m_s = \widetilde\bfg^0_sds$ in the first interval $[0,\frac{1}{N})$. The statement "\textit{The \rw $\{\widetilde\be_s\}_{s\geq 0}$ provides the solution of equation \eqref{EqDudleyDiff}}" can be given a precise meaning in the framework of non-standard analysis, and holds true, when correctly interpreted. This way of saying things is, in any case, useful (justified) and intuitive. 

\smallskip

\noindent \textbf{Notations.} We shall denote by $\widetilde\PP_{\be_0}$ the law of the \rw started from $\be_0$. Given two possibly infinite real numbers $a$ and $b$, we shall say that \textit{$a$ and $b$ are equal up to a negligeable quantity} if $\frac{a}{b}$ is infinitesimally close to $1$; we shall write $a\simeq b$. The notation $\textrm{Haar}(\cdot)$ will stand for a Haar measure on $SO_0(1,3)$.

\smallskip

$\bullet$ Formula \eqref{InterpretationOnePartMin} will hold true if we can prove it for any function $f$ on $\OO\VV_{t,\al}$ of the form
$$
f(m',\bfg') = {\bf 1}_A(m'){\bf 1}_G(\bfg'),
$$
for sufficiently small infinitesimal open sets $A\subset \VV_{t,\al}$ and $G\subset SO_0(1,3)$. We shall suppose, \wlg, that $A\times G$ is a (connected) neighbourhood of a given point $\be=(m,\bfg)\in\OO\VV_{t,\al}$. We shall associate to $\be$ the hyperplane $\VV_\be = \bigl\{m'\in\RR^{1,3}\,;\,m'\in m+\bigl(\bfg^0\bigr)^\perp\bigr\}$. To distinguish the Lebesgue measures induced by $q$ on $\VV_{t,\al}$ and $\VV_\be$, we shall denote them by $\LEB_3^{t,\al}$ and $\LEB_3^\be$ respectively.

If we let $A'$ be the set of points of $\VV_\be$ of the form $x+s{\bfg'}^0$, for $x\in A, s\in\RR$ and $\bfg'\in G$, the $\LEB_3^\be$-measure of $A'$ is equal to 
$$
\LEB_3^\be(V') \simeq q(\al^0,\bfg^0)\LEB_3^{t,\al}(A).
$$

Let now $M$ be an infinite integer and let run $M$ independent infinitesimal \rws started from $\be_0\in\OO\RR^{1,3}$. Write $N_{A\times G}$ and $N_{A'\times G}$ for the (random) numbers of trajectories of the \rw that hit $\OO V_{t,\al}$ and $\OO\VV_\be$ in $A\times G$ and $A'\times G$ respectively. If $A$ is small enough for $N \LEB_3^{t,\al}(A)$ to be infinitesimally close to $0$ and $M$ is large enough\footnote{Equal to an infinite integer depending on $N$ and $\LEB_3^{t,\al}(A)$.}, ($\PP_{\be_0}^{\otimes M}$-almost surely) 'almost all' the trajectories of the \rws hitting $A\times G$ will hit it in a time interval where $\bfg_s$ is constant. As the length of this time interval is much bigger than the time needed by any timelike path to go from $A$ to $A'$, the trajectories of the \rw will hit $A'\times G$ on the same time interval where they hit $A\times G$. As only a negligeable quantity of trajcetories hitting $A'\times G$ will not hit $A\times G$, we shall have on the one hand 
$$
N_{A\times G}\simeq N_{A'\times G}, \quad \widetilde\PP_{\be_0}^{\otimes M}-\textrm{almost surely}.
$$
As the strong law of large numbers ensures us that 
\begin{equation*}
\begin{split}
& N_{A\times G} \simeq M\times h^\al_{\be_0}(\be)\,\textrm{Haar}(G)\LEB_3^{t,\al}(A) \\
& N_{A'\times G} \simeq  M\times f(\be_0\,;\,\be)\,\textrm{Haar}(G)\LEB_3^\be(A') \simeq M\times q(\al^0,\bfg^0)f(\be_0\,;\,\be)\,\textrm{Haar}(G)\LEB_3^{t,\al}(A)
\end{split}
\end{equation*}
on the other hand, it follows that
$$
h^\al_{\be_0}(\be) \simeq q(\al^0,\bfg^0)f(\be_0\,;\,\be).
$$
Both quantities being standard reals, we actually have equality.

\bigskip
\bigskip

\noindent \textit{3.} We are now going to use equation \eqref{InterpretationOnePartMin} to give a proof of equation \eqref{FundamentalEq}. This will be done fixing a frame $\al\in SO_0(1,3)$ and proving that we have $\textstyle{L^*\Bigl(\frac{h^\la_{\be_0}}{\la_\al}\Bigr) = 0}$, where we have denoted by $\la_\al$ the function $\be=(m,\bfg)\mapsto q(\al^0,\bfg^0)$.

\smallskip

A frame $\al$ having been chosen, define the stopping times
$$
\forall\,r\in\RR,\quad S_r = \inf\{s>0\,;\,q(\al^0,m_s)>r\},
$$
and the process $\be^\al$ in the same way as the process $\be^\bfg$ has been defined above. It can be proved as above that the \rv $\be^\al_r$ has a smooth density $h_{\be_0}^\al(r,\cdot)$ \wrt $\Vol_{\OO\textrm{sp}(\al^1,\al^2,\al^3)}$, under $\PP_{\be_0}$; it is defined for $r\geq q(m_0,\al^0)$. The function $h_{\be_0}^\al(r,\cdot)$ is defined as equal to $0$ for $r<q(m_0,\al^0)$. Identifying $\bigl(r,(m',\bfg')\bigr)\in\RR\times\bigl(\OO\textrm{sp}(\al^1,\al^2,\al^3)\bigr)$ to the point $\bigl(r\al^0+m',\bfg'\bigr)$ of $\OO\RR^{1,3}$, the function $h_{\be_0}^\al$ will be seen as a function on $\OO\RR^{1,3}\backslash\{\be_0\}$. Three more notations will be needed: $L^\al$ will stand for the generator of the  $\OO\,\textrm{span}(\al^1,\al^2,\al^3)$-valued diffusion $\be^\al_.$, we shall write $D_mh^\al_{\be_0}$ for the partial differential of $h^\al_{\be_0}$ \wrt $m$(\footnote{The map ${D_mh^\al_{\be_0}}_{|\be}$ is for any $\be\in\OO\RR^{1,3}$ the linear form $\zeta\in\RR^{1,3}\mapsto \underset{\eta,0}{\lim} \frac{h^\al_{\be_0}(\be+\eta\zeta)-h^\al_{\be_0}(\be)}{\eta}$; this limit is denoted by $\bigl({D_mh^\al_{\be_0}}_{|\be}\bigr)(\zeta)$, or simply $(D_mh^\al_{\be_0})(\zeta)$.}) and use the notation $D_x$ to refer to the partial differentiation operation in the direction of $\textrm{span}(\al^1,\al^2,\al^3)$; last we shall decompose a vector $\bfg^0\in\HH$ as
$$
\bfg^0 = q(\bfg^0,\al^0)\al^0 + \sum_{i=1}^3\dot{x}_i\al^i.
$$
Note the relation
$$
\displaystyle{\bigl(D_m h^\al_{\be_0}\bigr) \left(\frac{\bfg^0}{q(\bfg^0,\al^0)}\right) = \partial_rh^\al_{\be_0} + \bigl(D_x h_{\be_0}^\al\bigr)\left(\frac{\dot{x}}{q(\bfg^0,\al^0)}\right)},
$$ 
which can be written
\begin{equation}
\label{EqDerivee}
-\left(D_x h^\al_{\be_0}\right)\left(\frac{\dot{x}}{q(\bfg^0,\al^0)}\right) = -H_0\left(\frac{h_{\be_0}^\al}{q(\bfg^0,\al^0)}\right) + \partial_rh_{\be_0}^\al.
\end{equation}

\smallskip

Recall that we write $\la_\al$ for $q(\al^0,\bfg^0)$. It can be proceeded like in the proof of the proposition/definition \ref{PropDefnOnePartDistrFunction} to show that $h^\al_{\be_0}(\cdot,\cdot)$ satisfies the heat equation
\begin{equation}
\label{HeatEqBis}
\partial_r\,h^\al_{\be_0} = \bigl(L^\al\bigr)^{*\al} h^\al_{\be_0},
\end{equation}
where 
$$
\bigl(L^\al\bigr)^{*\al} h^\al_{\be_0}= -\left(D_x h^\al_{\be_0}\right)\Bigl(\frac{\dot{x}}{\la_\al}\Bigr) + V^{*\al}\Bigl(\frac{h^\al_{\be_0}}{\la_\al}\Bigr) + \frac{1}{2}\bigl(V_iB^{ij}V_j\bigr)\Bigl(\frac{\la}{\la_\al}h^\al_{\be_0}\Bigr)
$$
and the operation ${}^{*\al}$ is the $\LL^2\bigl(\Vol_{\OO\,\textrm{sp}(\al^1,\al^2,\al^3)}\bigr)$-dual operation. Using equation \eqref{EqDerivee}, the heat equation \eqref{HeatEqBis} can be written
\begin{equation}
\label{LStarf=0}
-H_0\Bigl(\frac{h^\al_{\be_0}}{\la_\al}\Bigr) + V^{*\al} \Bigl(\frac{h^\al_{\be_0}}{\la_\al}\Bigr) + \frac{1}{2}\,\bigl(V_iB^{ij}V_j\bigr)\Bigl(\la\frac{h^\al_{\be_0}}{\la_\al}\Bigr) = 0.
\end{equation}
Note that since the vector field $V$ acts only on $SO_0(1,3)$ we have $V^{*\al} = V^*$; we have recalled above that $H_0^* = -H_0$. So, equation \eqref{LStarf=0} can take its final form: $L^*\Bigl(\frac{h^\al_{\be_0}}{\la_\al}\Bigr) = 0$, \textit{i.e.} $L^*f(\be_0\,;\,\cdot) = 0$.
\end{Dem}

This theorem/definition needs a few comments.

\smallskip

$\bullet$ Formula \eqref{FundamentalEq} is fundamental in the approach developped by Debbasch, Rivet and their co-workers. Their analysis of the situation entilery rests on a similar transport equation. Although it can be argued that since equation \eqref{InterpretationOnePartMin} implies that the one-particle distribution function determines the hitting distributions of the process at any times of any rest frame, a theorem of Blumenthal-Geetor and McKean ensures us that this function essentially determines the process, such a position should be taken with care. Indeed, the development of stochastic analysis has shown that one can gain much insight in the situation looking at the pathwise behaviour of processes rather that looking at analytic quantities such like their semi-group. We hope to illustrate this point throughout this article. In any case, theorem \ref{PropDefnOnePartDistrFunction} makes it clear that the fundamental quantity is not a hitting distribution $h^\al_{\be_0}$ but the one-particle distribution function; a fact which was not put forwards in the article \cite{ROUP01} of C. Barbachoux, F.Debbasch and J.P. Rivet.

\smallskip

$\bullet$ Equation \eqref{FundamentalEq} has a clear meaning from a Markov process point of view. It says that the measures $f(\be_0\,;\,\be)\Vol(d\be)$ on $\OO\RR^{1,3}$ are invariant for the $(V,\ff)$-diffusion. It is tempting to ask wether these measures and their possible renormalized limits as $\be_0$ goes to infinify are sufficient to describe the set of all invariant measures. For instance, it would be interesting, in the study of the R.O.U.P. in Minkowski space, to see if the strong recurrence of the process $\{\bfg_s\}_{s\geq 0}$ is sufficient to prove that the J\"uttner measure $ae^{-b\gamma}\Vol(d\be)$, alluded to above, is the only measure we obtain sending $\be_0$ at infinity, while imposing the limit measure to have mass in any open set\footnote{Other measures can be obtained if we do not impose this condition.}.
Even though a complete answer of the general question is out of reach at the moment, we shall come back in section \ref{PreBouquetFinal} to related matters in the general framework that we are going to present now.

\section{Relativistic diffusions in a Lorentzian manifold}
\label{DiffLorMan}

We shall now proceed to defining $(V,\ff)$-diffusions in a Lorentzian manifold. Let $(\MMM,q)$ denote a $(1+d)$-dimensional Lorentzian manifold, endowed with its Levi-Civita connection. As in Minkowski space, we shall construct the dynamics in a bigger space than $\MMM$. We shall first recall in section \ref{SubsubsectionDefnDiff}, \textbf{a)} how one can construct this space (the orthonormal frame bundle over $(\MMM,q)$) and the analogue of the above vector fields $H_0$ and $V_i$ before defining the class of $(V,\ff)$-diffusions in section \ref{SubsubsectionDefnDiff}, \textbf{b)}. We shall then see in section \ref{SubsubsectionSubDiff} that some situations give rise to a sub-diffusion in the (future-oriented) unit tangent bundle of $\MMM$. Several example will be discussed before returning in section \ref{PreBouquetFinal} to the study of $(V,\ff)$-diffusions. We shall define in this section the one-particle distribution function of the $(V,\ff)$-diffusion and prove its fundamental property. This result will shed some light on the structure of $L$-harmonic functions (section \ref{SubSubSectionHarmonicFunctions}) and will provide a simple proof of a general H-theorem (section \ref{BouquetFinal}). 

\medskip

\noindent \textbf{Hypothesis.} We shall suppose from now on that $(\MMM,q)$ is oriented and time-oriented.

\subsection{$(V,\ff)$-diffusions in $\OO\MMM$}
\label{SubsubsectionDefnDiff}

\paragraph{a) Geometrical objects in play.} 
Given some point $m\in\MMM$, it will be useful to consider an orthonormal basis $\{\bfg^0,...,\bfg^d\}$ of the tangent space $T_m\MMM$ to $\MMM$ at $m$ as an isometry from $\bigl(\RR^{1,3},q\bigr)$ to $\bigl(T_m\MMM,q\bigr)$(\footnote{The letter $q$ has here two different meanings.}); so, strictly speaking, $\bfg^i = \bfg(\varep^i)$.

The orthonormal frame bundle of $\MMM$ is just the collection 
$$
\OO\MMM = \bigl\{(m,\bfg)\,;\,m\in\MMM,\, \bfg \textrm{ an orthonormal basis of }T_m\MMM\bigr\}.
$$ 
We shall write  $\OO\mcU = \bigl\{(m,\bfg)\,;\,m\in \mcU,\, \bfg \textrm{ an orthonormal basis of }T_m\MMM\bigr\}$ for any subset $\mcU$ of $\MMM$. One defines the manifold structure of $\OO\MMM$ as follows. This structure being local, it suffices to define the structure of $\OO\mcU$ for any (small) domain $\mcU$ of $\MMM$; take it small enough to be the domain of a chart $x : \mcU \rightarrow \RR^{1+d}$. Applying Gram-Schmidt orthonormalisation procedure to the family of vectors $\bigl\{\partial_{x^i}\bigr\}_{i=0..3}$ in each tangent plane, on defines a section $\sigma : \mcU \rightarrow \OO\mcU$. The identification
$$
\bfi : \mcU\times O(1,3) \rightarrow \OO\mcU, \quad (m,g)\mapsto \bigl(m,\sigma (m)g\bigr)
$$
gives $\OO\mcU$ its differentiable structure (compatible with changes of charts)\footnote{Consult for instance chapter $10$ of the book \cite{MalliavinBook} of P. Malliavin.}. Note that $O(1,d)$ acts  on $\OO\MMM$ on the right: the action of $g'$ on $(m,\bfg)$ in the above chart $\bfi$ is
\begin{equation}
\label{ActionO}
(m,\bfg).g' = (m,\sigma (m)gg').
\end{equation}
Note that $\OO\MMM$ has several connected components. We shall be interested in dynamics leaving these components globally fixed. We choose to consider only one of them, specified by the requirement that $\bfg^0$ should be future-oriented and that the orientation of $\bfg$ should be direct (we have supposed the space oriented). The above action of the connected component of identity in $SO(1,d)$ preserves our connected component. We shall also denote it by $\OO\MMM$, as there will be no risk of confusion.

\smallskip

Action \eqref{ActionO} enables us to define vector fields on $\OO\MMM$:
$$
V_i\bigl((m,\bfg)\bigr) = \frac{d}{dt}_{\big|t=0}\bigl((m,\bfg).e^{tE_i}\bigr), \quad i=1..d.
$$
Last, we shall define the vector field $H_0$ as the infinitesimal generator of the geodesic flow on $\OO\MMM$. The dynamics $\bigl\{(m_s,\bfg_s)\bigr\}$ of this flow is described by asking that $\frac{dm_s}{ds} = \bfg^0_s$, and $\bfg_s$ should be transported parallely along the path $\{m_s\}$. One has for instance $H_0\bigl((m,\bfg)\bigr) = (\bfg^0,0)$ in Minkowski's flat spacetime, in accordance with the previous definition of $H_0$ given above.

\smallskip

\noindent \textbf{Notation.} We shall write $\be$ for a generic element of $\OO\MMM$.

\paragraph{b) $(V,\ff)$-diffusions.}
We are going to define $(V,\ff)$-diffusions following the same approach as in Minkowski space. We shall thus consider these diffusions as random perturbations of the flow of a differential equation in $\OO\MMM$ of the form
$$
d\be_s = H_0(\be_s)ds + V(\be_s)ds
$$
where $V$ is \textit{any} vector field on $\OO\MMM$. As in section \ref{ModelsInMin}, we shall not modelize the surrounding medium itself but just its action on the dynamics. This action will be given through the datum of  an $\OO\MMM$-valued previsible process $\{\ff_s\}_{s\geq 0}$ \st $\ff_s(\be_.) = \ff_s\bigl((m_.,\bfg_.)\bigr) = (m_s,f_s)$ for some orthonormal basis $f_s$ of $T_{m_s}\MMM$. Roughly speaking, it has the property that, \textit{when computed in the rest frame $\ff_s(\be_.)$, i.e. using its associated time}, the acceleration of $m_\cdot$ has a deterministic part and a random part which is Brownian in any spacelike direction belonging to $\textrm{span}\bigl(f^1_s(\be_.),f^2_s(\be_.),f^3_s(\be_.)\bigr)$.

Define $\{A_s\}_{s\geq 0}$ as the $d\times d$ random matrix process with coefficient $(i,j)\in [1,d]^2$ equal to $q(f^i_s,\bfg_s^j)$ at time $s$, and set
\begin{equation}
\label{DefnBeta1}
{\circ d\beta}_s = q(f_s^0,\bfg^0_s)^{\frac{1}{2}}\,A_s^{-1}{\circ dw_s}.
\end{equation}

\begin{defn}
\label{DefnDiff}
Define the $\RR^d$-valued process $\beta$ as above. A \textbf{$(V,\ff)$-diffusion in $(\MMM,q)$} is an $\OO\MMM$-valued process $\{\be_s\}_{s\geq 0} = \bigl\{(m_s,\bfg_s)\bigr\}_{s\geq 0}$ satysfying the stochastic differential equation
\begin{equation}
\label{EquationDynamics}
{\circ d}\be_s = H_0(\be_s)ds + V(\be_s)ds + V_i(\be_s)\,{\circ d\beta_s^i}.
\end{equation}
\end{defn}
If you do not feel comfortable with this stochastic differential equation, we shall give a step-by-step description of the dynamics in the next section. The remarks on probabilistic formalism and existence and uniqueness results made in section \ref{ModelsInMin}, \textbf{g)} apply here. Let us emphasize the interest that the above general definition might have for modelization. It provides a model of evolution of an object which has internal parameters (such as a spin) influencing the way it interacts with the surrounding medium, and whose value at some proper time depends on its past history. Challenging questions arise from this non-Markovianity of the model; yet, as only Markovian examples have been studied so far, we shall mainly explore this situation in the sequel.

\medskip

\noindent \textbf{Example: Franchi-Le Jan diffusion using co-ordinates.} This diffusion is the $(0,\be_.)$-diffusion, first defined in \cite{FranchiLeJan}. Note that the $(0,\be_.)$-diffusion is essentially the unique $(V,\ff)$-diffusion determined entirely by the geometric background $(\MMM,q)$. We asked
in section \ref{ModelsInMin}, \textbf{e)} which entity could give rise to the random excitement $V_i(\be_s)\,{\circ dw^i_s}$ that enters in the equations of motion of the Dudley(-Franchi-Le Jan)-diffusion process in the empty spacetime of Minkowski. This objection disappears when we consider the $(0,\be_.)$-diffusion in any spacetime $(\MMM,q)$ containing matter. It is in that case possible to add to the macroscopic description of matter given through the stress-energy-momentum (non-null) tensor a microscopic (quantum) description of matter from which randomness can be infered to come\footnote{Consult the article \cite{DowkerHensonSorkin} for results in this direction.}.

\smallskip

Equation \eqref{EquationDynamics} takes for this process the form
$$
{\circ d}\be_s = H_0(\be_s)ds + V_i(\be_s)\,{\circ dw_s^i}.
$$
To describe how we can write equation \eqref{EquationDynamics} using co-ordinates, note first that the data of local co-ordinates $x^i$ on $\MMM$ induces local co-ordinates on $T\MMM$: a vector $p\in T_m\MMM$ will be written $p = \displaystyle{\sum_{i=0..d}p^i\partial_{x^i}}$. Denoting then by $\Gamma : \RR^{1+d}\times\RR^{1+d} \rightarrow \RR^{1+d}$ the Christofel map associated with these co-ordinates, the dynamics of the Franchi-Le Jan diffusion takes the form 
\begin{equation}
\begin{split}
&      {\circ d}m_s = \bfg^0_sds, \\
& {\circ d}\bfg^0_s = -\Gamma(\bfg^0_s,\bfg^0_s)\,ds + \sum_{i=0..d}\bfg^i_s{\circ dw^i_s}, \\
& {\circ d}\bfg^j_s = -\Gamma(\bfg^0_s,\bfg^j_s)\,ds + \bfg^0_s{\circ dw^j_s}, \textrm{ for }j=1..d.
\end{split}
\end{equation}
These equations have to be written using the preceding co-ordinates. If one wishes to use Ito differentials, the system becomes
\begin{equation}
\label{EquationDynamicCoordIto}
\begin{split}
&      dm_s = \bfg^0_sds, \\
& d\bfg^0_s = \left(-\Gamma(\bfg^0_s,\bfg^0_s)+\frac{d}{2}\bfg^0_s\right)\,ds + \sum_{i=0..d}\bfg^i_s\, dw^i_s, \\
& d\bfg^j_s = \left(-\Gamma(\bfg^0_s,\bfg^j_s)+\frac{1}{2}\bfg^j_s\right)\,ds + \bfg^0_s\,dw^j_s, \textrm{ for }j=1..d.
\end{split}
\end{equation}
Remark that if we write $Q_m$ the matrix of the metric in these co-ordinates at point $m$, then the co-variance matrix of the martingale $\displaystyle{\sum_{i=0..d}\bfg^i_s\, dw^i_s}$ is equal to $\bfg^0_s(\bfg^0_s)^* - Q^{-1}_{m_s}$(\footnote{ We write here $\bfg^0$ for the vector of its co-ordinates in the basis $\{\partial_{x^i}\}_{i=0..d}$.}). The fact that it depends only on $m_s$ and $\bfg^0_s$ implies that the sub-process $\bigl\{(m_s,\bfg^0_s)\bigr\}_{s\geq 0}$ is itself a diffusion. The investigation of such situations is the object of the next section.

Note, en passant, that since we can read the matrix $Q_{m_s}$ on the co-variance of the martingale part of $\bfg^0_s$, it means that the \textit{local geometry} of $(\MMM,q)$ can be recovered from the pathwise study of the sub-process $\bigl\{(m_s,\bfg^0_s)\bigr\}_{s\geq 0}$. To determine what amount of information on the \textit{large scale structure} of the space $(\MMM,q)$ one can obtain from the pathwise study of this process or of the $(V,\ff)$-diffusion is a much harder task; we shall come back to it in section \ref{SubSubSectionHarmonicFunctions}.

\medskip

The heuristic explained in section \ref{ModelsInMin} and motivating the above definition of $(V,\ff)$-diffusions should make it clear that $(V,\ff)$-processes should be considered as models of diffusion in a homogeneous medium. Note yet that the input of a non-isotropic excitement in place of $dw$ in equations \eqref{DefnBeta1}, \eqref{EquationDynamics} would provide models of diffusions in a non-isotropic medium.

\medskip

\noindent \textbf{Hypothesis for the remainder of the article.} With in mind the diffusion of particles in a fluid, we shall suppose from now on that the flow of the vector field $V$ leaves each fiber of the projection $(m,\bfg)\in\OO\MMM \mapsto m\in\MMM$ stable.

\subsection{Sub-diffusions in $\HH\MMM$.} 
\label{SubsubsectionSubDiff}

As emphasized in section \ref{ModelsInMin}, \textbf{a)} in the framework of Minkowski spacetime, $(V,\ff)$-diffusions defined above can be considered as models of random motion of an infinitesimal rigid object in a relativistic medium. It might be interesting in some situations to define what could be the random motion of a point in such a medium. To investigate a physically motivated classical framework, we shall concentrate on Markovian processes.

As noted after definition \ref{DefnDiff}, the $(V,\ff)$-diffusions are not Markovian unless we choose a Markovian previsible process $\ff$:
$$
\ff_s(\be_.) = \ff(\be_s).
$$
This requirement is not sufficient yet to ensure that the sub-process $\bigl\{(m_s,\bfg^0_s)\bigr\}_{s\geq 0}$ of $\{\be_s\}_{s\geq 0}$ is itself a Markov process. We give in paragraph \textbf{a)} a simple condition which is proved to be sufficient in paragraph \textbf{b)}. Several examples are examined in paragraph \textbf{c)}.

\smallskip

Throughout this section, we shall suppose $\ff$ regular enough to have existence and strong uniqueness in the system \eqref{DefnBeta1}, \eqref{EquationDynamics}. We shall denote by\footnote{Recall that we have supposed $(\MMM,q)$ time-oriented.} 
$$
\HH\MMM = \bigl\{(m,\bfg^0)\in T\MMM\,;\,m\in\MMM,\,\bfg^0\in T_m\MMM \textrm{ future-oriented unit vector}\bigr\}
$$ 
the (future-oriented) unit sub-bundle of $T\MMM$. This space is the phase space of the set of $\mcC^1$ timelike paths in $(\MMM,q)$. The map 
$$
\pi : \OO\MMM \rightarrow \MMM
$$
will denote the projection $(m,\bfg) \mapsto m$, and $\widetilde\pi : \OO\MMM \rightarrow \HH\MMM$ the projection $(m,\bfg)\mapsto (m,\bfg^0)$.

\paragraph{a) A sufficient condition to a have a sub-diffusion in $\HH\MMM$.}
In addition to the hypothesis $\pi_* V = 0$ made above, we shall suppose that
\begin{itemize}
   \item there exists a vector field $\widehat{V}$ on $\HH\MMM$ \st $V$ is the horizontal lift of $\widehat{V}$ to $\OO\MMM$(\footnote{Denote by $\varphi_t(.)$ the flow of the vector field $V$ on $\OO\MMM$ and by $\{\widehat\varphi_t\}_{t\geq 0}$ the flow of $\widehat{V}$ on $\HH\MMM$. The above hypothesis means that the point $\varphi_t\bigl(m,(\bfg^0,\bfg^1,...,\bfg^d)\bigr)\in\OO_m\MMM$ is obtained by parallel transport of $(\bfg^1,...,\bfg^d)$ along the path $\bigl\{\widehat\varphi_s\bigl((m,\bfg^0)\bigr)\bigr\}_{s\leq t}$ in $\HH_m\MMM$.}).\end{itemize}

\smallskip

We shall begin our investigation with the particular case of the $(0,\be_.)$-diffusion of Franchi and Le Jan. Remember equation \eqref{ActionO} describing the action of $O(1,3)$ on $\OO\MMM$. This action induces a right action of $O(3)\subset O(1,3)$ on $\OO\MMM$, which amounts to rotate the vectors $\bfg^1, \bfg^2, \bfg^3$ in the Euclidean space they generate and leaves $\HH\MMM\subset\OO\MMM$ stable. Given the Brownian input $w$ in equations \eqref{DefnBeta1}, \eqref{EquationDynamics}, denote by 
\begin{equation}
\label{NotationSolution}
\be(s,\be_0\,;\,w)
\end{equation}
the (unique strong) solution started from $\be_0$. We have for any $g\in O(3)$
$$
\widetilde\pi\bigl(\be(s,\be_0g\,;\,w)\bigr) = \widetilde\pi\bigl(\be(s,\be_0\,;\,gw)\bigr).
$$
Since $gw$ is also a Brownian motion, the law of $\bigl\{\widetilde\pi\bigl(\be(s,\be_0\,;\,gw)\bigr)\bigr\}_{s\geq 0}$ does not depend on $g\in O(3)$, but only depends on $\widetilde\pi(\be_0)\in\HH\MMM$. The sub-process $\bigl\{\widetilde\pi\bigl(\be(s,\be_0\,;\,w)\bigr)\bigr\}_{s\geq 0}$ is thus a diffusion in $\HH\MMM$.

\paragraph{b) A step-by-step description of the dynamics.}
The general case is covered by the following theorem. 

\begin{thm}
\label{ThmDiffOnHM}
Suppose there exists a function $\ff^0 : \OO\MMM\rightarrow \HH\MMM$ \st 
\begin{equation}
\begin{split}
& \bullet\; \ff^0(\be)=\ff^0\bigl((m,\bfg)\bigr)\in\HH_m\MMM \textrm{ depends only on }(m,\bfg^0), \textrm{ and} \\
& \bullet\; \ff(\be_s)  = \Bigl(m_s,\bigl(\ff^0(\be_s),f^1(\be_s),...,f^d(\be_s)\bigr)\Bigr)
\end{split}
\end{equation}
for some functions $f^1,...,f^d$. Let $(m,\bfg^0)\in\HH\MMM$. Then, given any choice of $\be_0\in\OO\MMM$ \st $\widetilde\pi(\be_0) = (m,\bfg^0)$, the law of the $\HH\MMM$-valued process $\widetilde\pi\bigl(\be(.\,,\be_0\,;\,w)\bigr)$ depends only on $(m,p)$ and the function $\ff^0$, and not on the particular choice of $f^1,...,f^d$ and $\be_0$. The process $\widetilde\pi\bigl(\be(.,\be_0\,;\,w)\bigr)$ is a diffusion in $\HH\MMM$.
\end{thm}

We shall present a heuristic proof of this fact, the remaining work being just a matter of formalism. Equation \eqref{DefnBeta1} and \eqref{EquationDynamics} are the mathematical expression of the following heuristic dynamics explaining how on constructs $\be_{s+\delta s}$ from $\be_s$.

\begin{enumerate}
   \item Set $m_{s+\delta s} = m_s + \bfg^0_s\,ds$,
   \item then, set $\bfg^0_{s+\delta s} = \bfg^0_s + \delta\bfg^0_s + \widehat{V}_{(m_s,\bfg^0_s)}\delta s$. The increment $\delta \bfg^0_s$ is the only vector of $T_{\bfg^0_s}\bigl(\HH_m\MMM\bigr)$ \st its projection in $\textrm{span}\bigl(f^1(\be_s),...,f^d(\be_s)\bigr)$ parallelly to $\ff^0\bigl((m_s,\bfg^0_s)\bigr)$ is equal to the scaled Brownian increment $q\bigl(\ff^0(\be_s),\bfg^0_s\bigr)^{\frac{1}{2}}\,\displaystyle{\sum_{i=1}^d} f^i(\be_s)\,{\circ dw^i_s}$.
   \item Last, transport parallelly $\{\bfg^1_s,...,\bfg^d_s\}$ along the increment $\delta \bfg^0_s + \widehat{V}_{(m_s,\bfg^0_s)}$ of $\bfg^0_s$. 
\end{enumerate}
Examining this description of the dynamics, we see that any previsible orthonormal transform of the basis $\bigl\{f^1(\be_s),...,f^d(\be_s)\bigr\}$ will leave the law of the Brownian increment unchanged, so that the law of $\delta \bfg^0_s$ will also be left unchanged. Note also that the changing $\be_0\in\OO\MMM$ to another starting point with the same $\HH\MMM$-projection will only influence the dynamics of $\bfg^1_s,...,\bfg^d_s$. These remarks justify theorem \ref{ThmDiffOnHM}. To put this argument in a polished probabilistic form is a matter of formalism.

\paragraph{c) Examples.}
\begin{enumerate}
   \item \textbf{Dudley-Franchi-Le Jan diffusion in Minkowski spacetime (\cite{Dudley1}, \cite{FranchiLeJan}).} We have already given its description in \ref{ModelsInMin}, \textbf{f)}, $2$: $\{\bfg^0_s\}_{s\geq 0}$ is a \bm on $\HH$ and $m_s = m_0+\int_0^s \bfg^0_r\,dr$. The usual stochastic development procedure can be applied to this process to construct its $\HH\MMM$-version from its $\HH\RR^{1,3}$-version; see \cite{FranchiLeJan}, theorem $1$.
   
   \item \textbf{R.O.U.P. in Minkowski spacetime (\cite{ROUP0}).} This process is the $\HH\RR^{1,3}$-sub-process of the $(V,\textrm{Id})$-diffusion on $\OO\MMM$, where $V\bigl((m,\bfg)\bigr) = -\al\,\textrm{grad}(\ln \gamma)$, for some positive constant $\al$, and $\gamma = q(\varep^0,\bfg^0)$. In this flat spacetime with global co-ordinates $(t,x)$, the dynamics may be re-parametrized by the time $t$; the state space then becomes $\bigl\{(x,q)\in\RR^3\times\RR^3\bigr\}$, where $(t,x)$ are the co-ordinates of $m$ and $q$ is the $\textrm{span}(\varep^1,\varep^2,\varep^3)$-part of $\bfg^0$. With these notations, $\gamma = \gamma(q) = \sqrt{1+|q|^2_{\textrm{Eucl}}}$. The step-by-step description of the dynamics (or, more formally, the stochastic differential equation \eqref{EquationDynamics}) immediately yields the following stochastic differential equations for $(x_t,q_t)$, where $w$ is an $\RR^3$-Brownian motion:
\begin{equation}
\begin{split}
& dx_t = \frac{q_t}{\gamma(q_t)}dt, \\
& dq_t = -2\al\frac{q_t}{\gamma(q_t)}dt + {\circ dw_t};
\end{split}
\end{equation}
this is the original description of the R.O.U.P. up to some constants. 

Notice that the process $\{q_t\}_{t\geq 0}$ is a Kolmogorov diffusion in $\RR^3$. It has a unique invariant measure $\mu$, which is a probability and has density \wrt Lebesgue measure proportional to $e^{-4\,\al\,\gamma(q)}$. As we have $\underline\lim_{|q|,\infty} \bigl(|4\al\,\nabla\gamma|^2-4\al\triangle\gamma\bigr)(q) = 16\al^2 > 0$, a well known theorem ensures us that $\mu$ satisfies a Poincar\'e inequality. As is also well known\footnote{See for instance the book \cite{AneEtAl} of C. An\'e et al.}, this implies that the semi-group of the porcess $\{q_t\}_{t\geq 0}$ converges to equilibrium exponentially fast in $\LL^2(\mu)$, at least like $e^{-16(\al^2-\delta)\,t}$, for any $\delta>0$. This fact sheds some light on the numerical simulations made in section $4$ of the article \cite{ROUP0}.

\smallskip

\textbf{R.O.U.P. in an arbitrary inertial frame\footnote{Compare with the article \cite{ROUP01}.}.} It might be enlightening to write down the equation of the dynamics using the time $r$ and the $(\textrm{x,q})$-co-ordinates associated with any orthonormal frame $\bfg$ of $\RR^{1,3}$. We shall write $F(r,\textrm{q}_r)$ for the damping force in these co-ordinates; note that is depends on $r$ and $\textrm{p}_r$. Its precise expression is unimportant.

To take advantage of the description of the $\OO\RR^{1,3}$ process given in the above step by step description of the dynamics, and to take advantage of the irrelevance of the precise orthonormal frame $\{\bfg^1_s,\bfg^2_s,\bfg^3_s\}$ of $T_{\bfg^0_s}\HH$ we use in this construction, we chose to take as a basis of  $T_{\bfg^0_s}\HH$ the family 
$$
\bigl\{\bfg^1-(\bfg^1,\bfg^0_s)\bfg^0_s,\,\bfg^2-(\bfg^2,\bfg^0_s)\bfg^0_s,\,\bfg^3-(\bfg^3,\,\bfg^0_s)\bfg^0_s\bigr\}
$$
and write down the $\textrm{Vect}(\bfg^1,\bfg^2,\bfg^3)$-part of the increment of $d\bfg^0_s$ as
$$
(*) = \sum_{k=1..3}\left(\sum_{j=1..3}\bigl(\bfg^j-(\bfg^j,\,\bfg^0_s)\bfg^0_s,\bfg^k\bigr){\circ d\beta^j_s}\right)\bfg^k.
$$
In this expression, the matrix $A_s$ used to define $\beta$ has coefficient $(i,j)$ equal to $\left(\varep^i,\bfg^j-(\bfg^j,\bfg^0_s)\bfg^0_s\right)$; it depends only on $\bfg^0_s$. Switching from the description in terms of proper time $s$ to the evolution in terms of time $r$ results in multiplying $(*)$ by $\displaystyle{\Bigl(\frac{dr}{ds}\Bigr)^{\frac{1}{2}} =q(\bfg^0,\bfg^0_s)^{-\frac{1}{2}} = \gamma(\textrm{q}_r)^{-\frac{1}{2}}}$. This finally gives
\begin{equation}
\begin{split}
& d\textrm{x}_r = \frac{\textrm{q}_r}{\gamma(\textrm{q}_r)}dr, \\
& d\textrm{q}_r = F(r,\textrm{q}_r)\,dr + \gamma(\textrm{q}_r)^{-\frac{1}{2}}\,\sum_{k=1..3}\Bigl(\sum_{j=1..3}\bigl(\bfg^j-(\bfg^j,\,\bfg^0_r)\bfg^0_r,\bfg^k\bigr){\circ d\beta^j_r}\Bigr)\bfg^k.
\end{split}
\end{equation}
The vector $\bfg^0_r\in\HH$ is determined by $\textrm{q}_r$. No other choice of $\bfg^1_s, \bfg^2_s, \bfg^3_s$ would give something  fundamentally simpler. This complicated expression of the dynamics means nothing else than the inadequacy of the choice of co-ordinates to describe it.
 
   \item \textbf{R.O.U.P. in the spacially flat Robertson-Walker spacetime (\cite{ROUPCurved}).} This model on expanding universe is the product $\RR\times\RR^3$ equipped with a metric of the form $dt^2-a(t)^2dx^2$, where $a>0$. As in the preceding example, one can describe the trajectories of the R.O.U.P. using the absolute time $t$ and the state space $\bigl\{(x,q)\in\RR^3\times\RR^3\bigr\}$. We shall write $\gamma_t(q) = \sqrt{1+a(t)^2|q|^2_{\textrm{Eucl}}}$. The step-by-step scheme (or equation \eqref{EquationDynamics}) yields the equations of dynamics:
\begin{equation}
\begin{split}
& dx_t = \frac{q_t}{\gamma_t(q_t)}dt, \\
& dq_t = -2\al\,a(t)^2\frac{q_t}{\gamma_t(q_t)}dt + \frac{1}{a(t)}\,{\circ dw_t};
\end{split}
\end{equation}
the gradient part in $\HH_m\MMM$ gives rise to the term $\displaystyle{-2\al\,a(t)^2\frac{q_t}{\gamma_t(q_t)}dt}$, $w$ is an $\RR^3$-Brownian motion, and the $\displaystyle{\frac{1}{a(t)}\,{\circ dw_t}}$ term is the Brownian increment \textit{in the Euclidean space $\bigl(\RR^3,-a(t)^2dx^2\bigr)$}(\footnote{Compare the derivation of these equations with the approach of the article \cite{ROUPCurved} of F. Debbasch. Note that the dynamics is described in this paper not in $\HH\MMM\subset T\MMM$ but in $T^*\MMM$.}). It can be proved that this diffusion has an infinite lifetime.
   
   \item \textbf{Franchi-Le Jan diffusion in the spacially flat Robertson-Walker spacetime.} We shall use the notation $(m,\bfg^0)$ for a point of $\HH\MMM$ to describe the dynamics of this process. Using the canonical co-ordinates $(t,x)$ in $\RR\times\RR^3$ and denoting by $\bigl((t,x),(\dot{t},\dot{x})\bigr)$ the associated co-ordinates in $T\MMM$, we have seen in equation \eqref{EquationDynamicCoordIto} that the equations of the dynamics take the form\footnote{Consult for instance proposition $35$, p.$206$ of the book \cite{Oneill} by O'Neill for the computation of the Christoffel symbols in a warped product.}
\begin{equation}
\begin{split}
& dm_s = \bfg^0\,ds, \\
& d\dot{t}_s = \left(\frac{3}{2}\dot{t}_s-(aa')(t_s)\|\dot{x}_s\|^2_{\textrm{Eucl}}\right)\,ds + dM^{\dot{t}}_s, \\
& d\dot{x}_s = \left(\frac{3}{2}-2\frac{a(t_s)}{a'(t_s)}\dot{t}_s\right)\dot{x}_s\,ds + dM^{\dot{x}}_s,
\end{split}
\end{equation}
where the $\RR^4$-valued local martingale $M = \left(M^{\dot{t}},M^{\dot{x}}\right)$ has co-variance
$$
\begin{pmatrix} \dot{t}_s^2-1 & \dot{t}_s\dot{x}_s^* \\ \dot{t}_s\dot{x}_s & \dot{x}_s\dot{x}_s^* + a^{-2}(t_s)\textrm{Id}_3 \end{pmatrix}.
$$
It has been shown by J. Angst that this diffusion has an infinite lifetime. Remark that the $\RR^2$-valued sub-process $\bigl\{(t_s,\dot{t}_s)\bigr\}_{s\geq 0}$ is a diffusion. This kind of decomposition of the diffusion into smaller dimensional diffusions has been the key of the previous investigations in Schwarzschild and G\"odel's spacetimes. See \cite{FranchiLeJan} and the article \cite{Franchi} of J. Franchi.


\end{enumerate}

\subsection{One-particle distribution function}
\label{PreBouquetFinal}

The aim of this section is to clarify the so-called notion of one-article distribution function in the general framework of Markovian $(V,\ff)$-diffusions on any Lorentzian manifold. We shall define it properly in section \ref{SubSubSectionOnePartFunction} and prove in theorem \ref{LStarOnePartFunction} that it satisfies a remarkable equation. This theorem will justify the analytic approach developped by F. Debbasch and his co-workers, as exposed in \cite{ROUPUnifying1} or \cite{ROUPUnifying2} and the references cited therein. The relevance of this notion in the study of $L$-harmonic functions will be described in section \ref{SubSubSectionHarmonicFunctions}.

\smallskip

The approach to one-particle distribution functions developped here is similar in spirit to the physical approach exposed in the article \cite{DRVanLeeuwen}, in a physical/mathematical style\footnote{Consult also the article \cite{IsraelBoltzmann} of W. Israel for a similar point of view.}. It should be noted yet that only the special relativistic situation is investigated in this article, whereas we deal below with the general relativistic case.

\smallskip

In order to ease the understanding of the situation, we shall make a hypothesis on the \textit{global geometry} of the space $(\MMM,q)$. We shall suppose the spacetime$(\MMM,q)$ \textbf{\textit{strongly causal}}: every point of $\MMM$ has arbitrary small (connected) neighbourhoods which no non-spacelike paths intersect more than once. This is a mild global assumption on the geometry of the space, satisfied by most of the models of physical spacetimes. This excludes, yet, pathological spaces where closed timelike paths exist, like G\"odel's spacetime.

We shall also use the following \textit{local} property, shared by all Lorentzian open manifold. Any point has an open (relatively compact connected) neighbourhood on which a time function is defined. By time function we mean a smooth function whose level sets are spacelike hypersurfaces\footnote{We shall construct these neighbourhoods in the beginning of the proof of proposition/definition \ref{DefnOnePartDistr}.}. We shall denote by $\mcU_m$ such a neighbourhood associated to a point $m\in\MMM$. We shall suppose $\mcU_m$ small enough to have the property that no non-spacelike paths intersect it more than once. As a consequence, it will enjoy the following extra property. Any timelike path in $(\MMM,q)$ will hit any spacelike hypersurface of $\mcU_m$ at most once. This property will be the main ingredient used to define of the one-particle distribution of $(V,\ff)$-diffusions.

\subsubsection{One-particle distribution function.}
\label{SubSubSectionOnePartFunction} 
The initial point $\be_0$ of the $(V,\ff)$-diffusion will be fixed throughout this paragraph. Given a point $\be=(m,\bfg)\in\OO\MMM$, different from $\be_0$, define the collection 
$$
\mcalV_\be = \bigl\{\VV\,;\,\textrm{ spacelike hypersurfaces of }\MMM \textrm{ contained in }\mcU_m \textrm{ and \st } m\in\VV \textrm{ and }T_m\VV = \bigl(\bfg^0\bigr)^{\perp}\bigr\}.
$$
Associate to any $\VV\in\mcalV_\be$ the hitting time
$$
H = \inf\{s\geq 0\,;\,m_s\in\VV\}.
$$
Given any point $m'$ in $\MMM$ we shall denote by $\Vol_{m'}(d\bfg)$ the Haar measure on $\OO_{m'}\MMM$, normalized in such a way that its projection on $\HH_{m'}\MMM$ is the Riemannian volume element induced by $q$. Recall the definition of $\OO\VV = \bigl\{(\widehat m,\widehat \bfg)\in\OO\MMM\,;\,\widehat m\in\VV, \widehat \bfg\in T_m\MMM\bigr\}$. Let us insist on the fact that even if $\VV$ is a sub-manifold of $\MMM$, the element $\bfg'$ of a point $(\widehat m,\widehat \bfg)\in\OO\VV$ is \textit{not} an orthonormal basis of $T_{\widehat m}\VV$, but an orthonormal basis of $T_{\widehat m}\MMM$. We shall write $\widehat\be = (\widehat m,\widehat\bfg)$ for a generic element of $\OO\VV$ and shall denote by $\sigma_\VV(d\widehat{m})$ the volume element induced by $q$ on $\VV$. With these notations, we shall endow the bundle $\OO\VV$ with the measure
$$
\Vol_{\OO\VV}(d\widehat\be) = \Vol_{\widehat{m}}(d\widehat\bfg)\otimes\sigma_\VV(d\widehat{m}).
$$
Recall that the point $\be=(m,\bfg)\in\OO\MMM$ has been fixed above.

\begin{PropDefn}
\label{DefnOnePartDistr}
Let $\VV\in\mcalV_\be$.
\begin{enumerate}
   \item The \rv $\be_H{\bf 1}_{H<\infty}$ has a smooth density $f_\VV(\be_0\,;\,\widehat\be)$ \wrt the measure $\Vol_{\OO\VV}(d\widehat\be)$ on $\OO\VV$.
   \item We have $f_{\VV'}(\be_0\,;\,\be) = f_\VV(\be_0\,;\,\be)$ for any other $\VV'$ in $\mcalV_\be$.
\end{enumerate}

\smallskip

\noindent So this quantity $ f_\VV(\be_0\,;\,\be)$ is independent from $\VV\in\mcalV_\be$; call it the value at point $\be$ of the \textbf{one-particle distribution function of the $(V,\ff)$-diffusion started from $\be_0$}. We shall denote it by $f(\be_0\,;\,\be)$; it is defined for $\be\neq\be_0$.
\end{PropDefn}

As is clear from its definition, this function takes the same value on points with the same $\HH\MMM$-projection. Given any point $\be_0\in\OO\MMM$, we shall adopt the usual conventions and shall denote by 
$$
I^+(\be_0) = I^+\bigl((m_0,\bfg_0)\bigr) = \bigl\{(\gamma(1),\bfg')\in\OO\MMM\,;\,\gamma \textrm{ future-oriented timelike path}, \gamma(0)=m_0,\,\bfg'\in\OO_{\gamma(1)}\MMM\bigr\}
$$
the chronological future of $\be_0$. This is an open set of $\OO\MMM$. It comes from the support theorem of Stroock and Varadhan that $f(\be_0\,;\,\cdot)$ is positive in $I^+(\be_0)$ and null outside the closure of $I^+(\be_0)$(\footnote{In Minkowski spacetime, this result comes from proposition $8$ in the article \cite{BailleulPoisson} of I. Bailleul. A similar proof can be given in the general framework of Markovian $(V,\ff)$-diffusions on any $\OO\MMM$.}).

\smallskip

We describe here the proof of proposition/definition \ref{DefnOnePartDistr} without technicalities. 

\begin{figure}[h!]
\label{FigureOnePartDistrFunction}
\begin{center}
\input{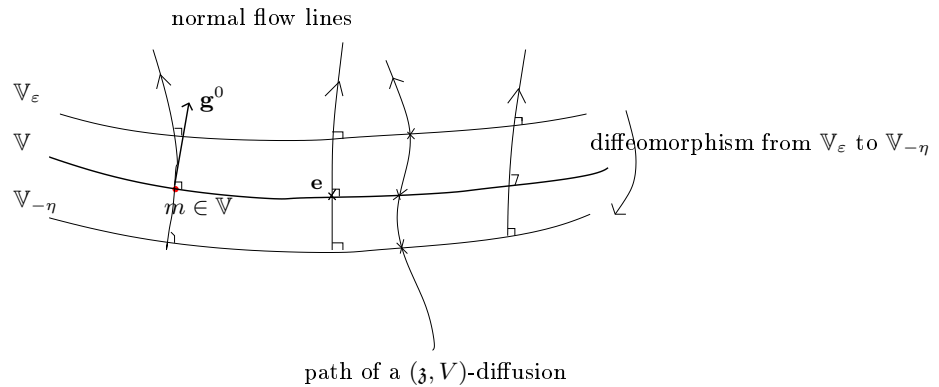} 
\end{center}
\caption{Constructing the one-particle distribution function.}
\end{figure}
We use the same idea as in section \ref{SectionOnePart} where a family of constant time hyperplanes was used to re-parametrize the process. These global objects will be here replaced by local ones: the normal variation $\{\VV_\varep\}_{\varep\in (-\eta,\eta)}$ of the spacelike hypersurface $\VV$. Their local definition is illustrated in figure $3$. Suppose that $\be_0$ belongs to $\VV_{-\eta}$; the timelike path $\{\be_s\}_{s\geq 0}$ will then hit each $\VV_\varep$ once, at increasing $\varep$. So we can use $\varep$ as a time parameter in place of $s$. Using the flow of the normal variation, we can consider the re-parametrized path as a hypoelliptic diffusion in $\VV$. We shall then get the conlcusion from H\"ormander's theorem on hypoellipticity.

\smallskip

The second point of proposition/definition \ref{DefnOnePartDistr} is established using the property of the sets $\mcU_m$ mentionned in the introduction. Indeed, suppose that $\VV'$ and $\VV$ not only have the same tangent space at $m$ but are equal in a neighbourhood $U$ of $m$ in $\VV$. Then, since $\VV'$ is a spacelike hypersurface of $\MMM$ contained in $\mcU_m$ and since no spacelike path of $(\MMM,q)$ can hit $\VV'$ or $\VV$ more than once, the set of trajectories of the $(V,\ff)$-process hitting $\OO\VV'$ in $\OO U$ is the same as the set of trajectories hitting $\OO\VV$ in $\OO U$. So the densities of $\be_{H_{\VV'}}{\bf 1}_{H_{\VV'}<\infty}$ and $\be_{H_{\VV}}{\bf 1}_{H_{\VV}<\infty}$ are under $\PP_{\be_0}$ equal on $\OO U$, \textit{i.e.} $f_{\VV'}(\be_0\,;\,\widehat\be) = f_\VV(\be_0\,;\,\widehat\be)$ for any $\widehat\be\in U$. Shrinking $U$ to $\{m\}$ formally gives the second point of proposition/definition \ref{DefnOnePartDistr}.

\begin{Dem}
\textit{1.} \textbf{Normal variation of a spacelike hypersurface.} Let $\VV\in\mcalV_\be$. For $\widehat m\in\VV$ and $\varep$ small enough, define $\phi_\varep(\widehat m)$ as the position at time $\varep$ of the geodesic started from $\widehat m$ leaving $\VV$ orthogonally, in the future direction. Then there exists (as a consequence of the local inversion theorem) a positive constant $\eta$ and an open set $\mcU\subset\MMM$ \st the map $\phi : (-\eta,\eta)\times\VV \rightarrow\mcU$, $(\varep,\widehat m)\mapsto \phi_\varep(\widehat m)$ is a diffeomorphism. It has the following properties, where we write $\VV_\varep$ for $\phi_\varep(\VV)$.
\begin{itemize}
   \item $\phi_0(\widehat m) = \widehat m$,
   \item $\partial_\varep\phi_\varep(\widehat m)\in\HH_{\phi_\varep(\widehat m)}\MMM$, and
   \item $\partial_\varep\phi_\varep(\widehat m)$ is orthogonal to $T_{\phi_\varep(\widehat m)}\VV_\varep$.
\end{itemize}
The family of spacelike hypersurfaces $\{\VV_\varep\}_{\varep\in(-\eta,\eta)}$ is called the \textbf{\textit{normal variation of}} $\VV$. The open set $\mcU$ has the function $\varep$ as a time function\footnote{These sets $\mcU$ are those used in the introduction to construct the sets $\mcU_m$.}. We shall suppose \wlg that $\mcU$ is diffeomorhic to an open set of $\RR^{1+d}$. The diffeomorphism $\phi$ can be extended to $(-\eta,\eta)\times\OO\VV \rightarrow\OO\mcU$. To that end, given $\varep\in(-\eta,\eta)$ transport parallelly $\bfg\in\OO_{\widehat m}\MMM$ along the path $\bigl\{\phi_t(\widehat m)\bigr\}_{t\in [0,\varep]}$; write $T^\phi_{\varep\leftarrow 0}\bfg$ for the element of $T_{\phi_\varep(\widehat m)}\MMM$ obtained that way. The map 
$$
\bigl(\varep,(\widehat m,\bfg)\bigr)\in(-\eta,\eta)\times\OO\VV \rightarrow \bigl(\phi_\varep(\widehat m),T^\phi_{\varep\leftarrow 0}\bfg\bigr)\in\OO\mcU
$$
is easily seen to be a diffeomorphism extending $\phi$. We shall still denote it by $\phi$.

\smallskip

\textbf{Notations.} Given a point $m\in\VV_\varep$, we shall denote by $\varpi(m)$ the future unit timelike vector orthogonal to $T_m\VV_\varep$. We can extend this vector fields $\varpi$ on $\mcU\subset\MMM$ to a vector field on $\OO\mcU$ lifting it horizontally; we shall still denote it by $\varpi$. In addition to this vector field $\varpi$ on $\OO\mcU$ we shall need some more notations.
\begin{itemize}
   \item $\overline\gamma := q(\varpi(\be),\bfg^0)$ will be a function of $\be=(m,\bfg)\in\OO\MMM$.
   \item The $^{*\OO\VV_\varep}$-operation will stand for taking the $\LL^2(\Vol_{\OO\VV_\varep})$-dual.
   \item Last, $H_\varep$ will denote the hitting time of $\OO\VV_\varep\subset\OO\MMM$.
\end{itemize}

Given a point $\be\in\OO\MMM$, we shall denote by $I^-(\be)$ its \textit{\textbf{timelike past}}:
$$
I^-(\be) = I^-\bigl((m,\bfg)\bigr) = \bigl\{(\gamma(1),\bfg')\in\OO\MMM\,;\,\gamma \textrm{ past-oriented timelike path}, \gamma(0)=m,\,\bfg'\in\OO_{\gamma(1)}\MMM\bigr\}.
$$
The timelike past of a set will be the union of the timelike past of its elements.

\medskip

\textbf{\textit{a)}} We shall suppose first that $\be_0$ belongs to $\OO\mcU$. If $\be_0=(m_0,\bfg_0)$ does not belong to the (closure of the) timelike past of $\VV$, then no timelike path started from $m_0$ can ever hit $\VV$, so the function $f(\be_0\,;\,\cdot)$ is null in a neighbourhood of $\OO\VV$. As we are interested in what happens near $\OO\VV$, we shall make the hypothesis that $\be_0$ belongs to the timelike past of $\VV$. We shall suppose, \wlg, that $\be_0\in\VV_{-\eta}$; it will be fixed throughout this paragraph.

\smallskip

As the hitting times $H_\varep$ will be $\PP_{\be_0}$-almost surely finite under the preceding hypothesis, we can consider the re-parametrized process $\{\be_{H_\varep}\}_{\varep\in(-\eta,\eta)}$; it has generator $\overline\gamma^{-1}\,L$. We shall decompose this operator under the form
\begin{equation}
\label{DecompositionGenerator}
\forall\,\be=\phi_\varep(\widehat\be),\quad \frac{Lf}{\overline\gamma}(\be) = (\varpi f)(\be) + \widehat L(f\circ \phi_\varep)\,(\widehat\be) = (\varpi f)(\be) + \bigl(\overline L f\bigr)(\be),
\end{equation}
where $\widehat L$ is a (smooth) second order differential operator on $\OO\VV$, and where, as a consequence, $\overline L$ acts only on $\OO\VV_\varep$. Now, define the $\OO\VV$-valued process $\bigl\{\widehat\be_\varep\bigr\}_{\varep\in(-\eta,\eta)} := \bigl\{\phi_\varep^{-1}(\be_{H_\varep})\bigr\}_{\varep\in(-\eta,\eta)}$ and denote by $\widehat\ell_\varep$ its time-dependent generator. The vector fields $V$ and $V_i$ acting only on the fibers of the projection $\OO\MMM\rightarrow\MMM$, it is easily seen that $\widehat\ell_\varep$ is a hypoelliptic operator, so the \rv $\widehat\be_\varep$ has for any $\varep\in (-\eta,\eta)$ a smooth density \wrt to $\Vol_{\OO\VV}$. It follows that $\be_H = \widehat\be_0$ also has a smooth density \wrt $\Vol_{\OO\VV}$.

\smallskip

\textbf{\textit{b)}} To deal with the general case where $\be_0$ does not belong to $\OO\mcU$, we can suppose \wlg that $\VV$ is a subset of a spacelike hypersurface $\VV'$ \st the analysis of point \textbf{\textit{a)}} applies and \st any timelike path hitting $\OO\mcU$ hits $\VV'_{-\eta'}$ before. Then, denoting by $h(\be_0\,;\,\widehat\be)\Vol_{\OO\VV'}(d\widehat\be)$ the smooth hitting distribution of $\OO\VV'$ by the process $\be_.$ under $\PP_{\be_0}$, we have
$$
\forall\,A\subset \OO\VV,\quad \PP_{\be_0}\bigl(\be_H\in A,\, H<\infty\bigr) = \int_A \left(\int h(\be_0\,;\,\widehat\be)f(\widehat\be\,;\,\be')\Vol_{\OO\VV'_{-\eta'}}(d\widehat\be)\right)\Vol_{\OO\VV}(d\be');
$$
from which we conclude that the \rv $\be_H{\bf 1}_{H<\infty}$ has under $\PP_{\be_0}$ a smooth density \wrt $\Vol_{\OO\VV}$, equal to $f_\VV(\be_0\,;\,\be) = \int h(\be_0\,;\,\widehat\be)f(\widehat\be\,;\,\be')\Vol_{\OO\VV'_{-\eta'}}(d\widehat\be)$.

\bigskip

\textit{2.} The formal proof of this point proceeds using a slightly different point of view than the heuristic described before the beginning of the proof of proposition/definition \ref{DefnOnePartDistr}. Fix $\be=(m,\bfg)\in\OO\MMM$ and let $\eta_0>0$ be smaller than the radius of definition of the (Lorentzian) exponential map $\exp_m : T_m\MMM\rightarrow \MMM$, and small enough for the geodesic ball of radius $\eta_0$ to be included in $\mcU_m$. Given $\eta<\eta_0$, denote by $A_\eta$ the hypersurface of $\MMM$ defined as 
$$
A_\eta := \{\exp_m(sT)\,;\,|s|<\eta,\,T\in\bigl(\bfg^0\bigr)^\perp\}.
$$
For $\eta_0$ small enough, the hypersurface $A_{\eta_0}$ will be spacelike; pick such an $\eta_0$. Denote also by $B_\eta$ the set of points of $\MMM$ of the form $\exp_{m'}(sU)$ for $m'\in A_\eta, |s|<\eta^2$ and $U\in T_{m'}\MMM$.

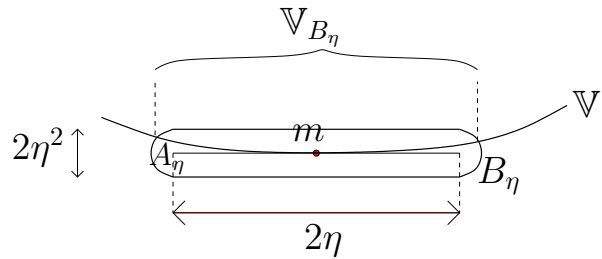
\begin{figure}[h!]
\label{DemoOnePartFunction}
\begin{center}
\input{DemoOnePartFunction.pstex_t} 
\end{center}
\end{figure}

The set $B_\eta$ has two important properties. We use the notation $\VV$ for any spacelike hypersurface belonging to $\mcalV_\be$. Recall that $\sigma_\VV$ stands for the volume element induced by $q$ on $\VV$.
\begin{equation}
\label{PremiereAsymptotique}
\frac{\Vol(B_\eta)}{\eta^{d+2}}\underset{\eta,0^+}{\longrightarrow} c_d
\end{equation}
If we write $\VV_{B_\eta}$ for the intersection of $\VV$ with the chronological past and future of $B_\eta$ in $\mcU_m$, we have
\begin{equation}
\label{DeuxiemeAsymptotique}
\frac{\sigma_\VV(\VV_{B_\eta})}{\eta^d}\underset{\eta,0^+}{\longrightarrow} c_d.
\end{equation}
The constant $c_d$ appearing above is the Euclidean volume of the unit ball of $\RR^d$. Given now any hypersurface $\VV\in\mcalV_\be$, $0<\eta<\eta_0$ and a positive integer $N$, run $N$ independent $(V,\ff)$-diffusions started from $\be_0$. We shall write $\be^{(i)}_{H^{(i)}}$ for the random position of the $i^{\textrm{th}}$ diffusion stopped at the random time $H^{(i)}$ where it hits $\OO\VV$ (provided this time is finite). Associate to a given real valued Lipschitz function $\varphi$ on $\OO\mcU_m$ the \rv
\begin{equation}
\label{TroisiemeAsymptotique}
F_N(\eta) := \sum_{i=1..N} \varphi\bigl(\be^{(i)}_{H^{(i)}}\bigr){\bf 1}_{H^{(i)}<\infty}{\bf 1}_{\be^{(i)}_{H^{(i)}}\in\OO\VV_{B_\eta}}.
\end{equation}
The almost sure following limit is a consequence of the strong law of large numbers:
$$
\underset{N+\infty}{\lim} \frac{F_N(\eta)}{N} = \int_{\OO\VV_{B_\eta}} \varphi(\widehat\be)\,f_\VV(\be_0\,;\,\widehat\be)\Vol_{\OO\VV}(d\widehat\be).
$$
If we now let $H_{B_\eta}^{(i)}$ be the hitting time of the set $B_\eta$ by the $i^{\textrm{th}}$ $(V,\ff)$-diffusion, set
$$
G_N(\eta) := \sum_{i=1..N} \varphi\Bigl(\be^{(i)}_{H_{B_\eta}^{(i)}}\Bigr){\bf 1}_{H_{B_\eta}^{(i)}<\infty}.
$$
Since $\varphi$ is Lipschitz and $B_\eta$ has a 'height' of order $\eta^2$, we have 
$$
\left|\varphi\Bigl(\be^{(i)}_{H^{(i)}}\Bigr){\bf 1}_{H^{(i)}<\infty}{\bf 1}_{\be^{(i)}_{H^{(i)}}\in\OO\VV_{B_\eta}} - \varphi\Bigl(\be^{(i)}_{H_{B_\eta}^{(i)}}\Bigr){\bf 1}_{H_{B_\eta}^{(i)}<\infty}\right| \leq C \eta^2
$$
for some positive constant $C$; it follows that
$$
\left|\underset{N,\infty}{\underline{\overline{\lim}}} \left(\frac{G_N(\eta)}{N}-\frac{F_N(\eta)}{N}\right)\right| \leq C\eta^2.
$$
Together with equations \eqref{PremiereAsymptotique}, \eqref{DeuxiemeAsymptotique}, \eqref{TroisiemeAsymptotique}, this equation gives us the existence and the value of the limit
$$
\underset{\eta,0}{\lim}\,\left(\Vol(B_\eta)^{\frac{-d}{d+2}}\,\underset{N,\infty}{\lim}\frac{G_N(\eta)}{N}\right) = \int_{\OO_m\MMM} \varphi(m,\widehat\bfg)\,f_\VV\bigl(\be_0\,;\,(m,\widehat\bfg)\bigr)\Vol_m(d\widehat\bfg).
$$
The left hand side being independent of $\VV$, the functionnal of the Lipschitz function $\varphi$ defined by the right hand side is also independent of $\VV$. The class of Lipschitz functions if rich enough to conclude from that fact that the measure $f_\VV(\be_0\,;\,\cdot)\Vol_m(\cdot)$ is independent of $\VV\in\mcalV_\be$, which implies that $f_\VV(\be_0\,;\,\be)$ itself is  independent of $\VV\in\mcalV_\be$.
\end{Dem}

To state the next proposition on $f(\be_0\,;\,\cdot)$ we shall write $\VV$ for a spacelike hypersurface of $\MMM$ and shall denote by $H$ the hitting time of $\OO\VV$. Given a point $\widehat\be=\bigl(\widehat m,\widehat\bfg\bigr)\in\OO\VV$, we shall denote by $\varpi_\VV(\widehat\be)$ the future unit timelike vector orthogonal to $T_{\widehat m}\VV$ (in accordance with the previous notation). This fundamental proposition extends the second point of theorem \ref{PropDefnOnePartDistrFunction} to the general framework adopted in this section.

\begin{prop}
\label{InterpretationOnePart}
Let $\be_0$ be a point of $\OO\MMM$ not belonging to $\OO\VV$. We have 
\begin{equation}
\label{EqInterpretationOnePartFunction}
\EE_{\be_0}\bigl[f(\be_H){\bf 1}_{H<\infty}\bigr] = \int_{\OO\VV} f(\widehat\be)\,q\bigl(\widehat\bfg^0,\varpi_\VV(\widehat\be)\bigr)f(\be_0\,;\,\widehat\be)\,\Vol_{\OO\VV}(d\widehat\be)
\end{equation}
for any bounded function $f$ on $\OO\VV$.
\end{prop}

\begin{DemOuv}
We shall use the notation $f_\VV(\be_0\,;\,\cdot)$ to denote the (smooth) density of the law of the \rv $\be_H {\bf 1}_{H<\infty}$ under $\PP_{\be_0}$, \wrt $\Vol_{\OO\VV}$. Given a point $\be\in\OO\VV$, we are going to prove that
\begin{equation}
\label{FundRel}
f_\VV(\be_0,\;\,\be) = q\bigl(\varpi_\VV(\be),\bfg^0\bigr)f(\be_0\,;\,\be).
\end{equation}
This point $\be=(m,\bfg)$ is now fixed. We shall denote by $\WW$ a hypersurface of $\mcalV_\be$; we have seen in proposition/definition \ref{DefnOnePartDistr} that $f(\be_0\,;\,\be) = f_\WW(\be_0\,;\,\be)$.

\smallskip

\textbf{Idea of the proof.} The idea of the proof is simple and illustrated in figure $4$. Pick a positive integer $N$; it will be sent to infinity at the end of the proof. Let $\bigl\{\VV_\varep\bigr\}_{\varep\in (-\eta,\eta)}$ be the normal variation of $\VV$. The positive real $\eta$ is chosen in such a way that any timelike geodesic started from $\VV_{-\eta}$, of length $\geq\frac{1}{N}$, hits $\VV_\eta$. It implicitly depends on $N$; we choose it as a decreasing function of $N$ converging to $0$ as $N$ increases to infinity. We shall write $\VV_{-\eta}'$ for the set of points of $\VV_{-\eta}$ from which any future-oriented timelike path hits $\VV$. We can suppose \wlg that the hypersurface $\WW$ is included in $\mcU = \displaystyle{\bigcup_{\varep=(-\eta..\eta)}\VV_\varep}$.

\begin{figure}[h!]
\label{PropFundOnePartFction}
\begin{center}
\input{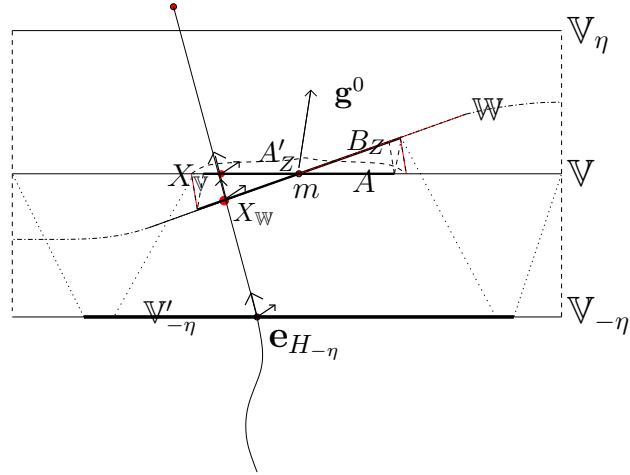} 
\caption{Proof of proposition \ref{InterpretationOnePart}.}
\end{center}
\end{figure}
Given a point $\widetilde\be = (\widetilde m,\widetilde\bfg)\in\OO\VV_{-\eta}'$, we shall write $m_\VV(\widetilde\be)$ for the intersection of the future-oriented geodesic $\gamma_{\widetilde\be}$ started from $\widetilde m$ in the direction $\widetilde\bfg^0$ and by $\bfg_\VV(\widetilde\be)$ the image at the point $m_\VV(\widetilde\be)$ of $\widetilde\bfg$ by parallel transport along $\gamma_{\widetilde\be}$. We set
$$
X_\VV(\widetilde\be) := \bigl(m_\VV(\widetilde\be),\bfg_\VV(\widetilde\be)\bigr).
$$
The point $X_\WW(\widetilde\be)$ is defined similarly using $\WW$ in place of $\VV$. Let us denote by $H_{-\eta}'$ the hitting time of $\OO\VV_{-\eta}'$ and set 
$$
X_\VV := X_\VV(\be_{H_{-\eta}'}){\bf 1}_{H_{-\eta}'<\infty} \in\OO\VV \textrm{   and   } X_\WW := X_\WW(\be_{H_{-\eta}'}){\bf 1}_{H_{-\eta}'<\infty}\in\OO\WW.
$$
We are going to see that these random points have smooth densities at $\be$ which satisfy equation \eqref{FundRel}; they both depend on $N$. Equation \eqref{FundRel} itself will be obtained as a limit, sending $N$ to infinity.

\medskip

\textbf{Proof.} Given a (small) open neighbourhood $\underline{A}$ of $m$ in $\VV$ define $\underline{B}\subset\WW$ as the intersection of $ \WW$ with the chronological past and future of $\underline A$ in $\mcU$. Pick $\underline A$ and $\eta$ small enough in such a way that any timelike path hitting $\underline A$ or $\underline B$ hits $\VV_{-\eta}'$ before. All open sets $A$ used hereafter will implicitly be supposed to be included in this fixed $\underline A$. 

The following lemma is proved noting that the maps $X_\VV : X_\VV^{-1}(\OO\underline A) \rightarrow \OO\underline A$ and $X_\WW : X_\WW^{-1}(\OO\underline B) \rightarrow \OO\underline B$ are well defined smooth diffeomorphisms; as such, the push forwards of any smooth measure on $\OO\VV'_{-\eta}$ by these maps are smooth measures on $\OO\underline A$ and $\OO\underline B$ respectively.

\smallskip

\begin{lem}
The laws of the \rvs $X_\VV{\bf 1}_{X_\VV\in\OO \underline A}$ and $X_\WW{\bf 1}_{X_\WW\in\OO \underline B}$ under $\PP_{\be_0}$ have smooth densities \wrt ${\bf 1}_{\OO\underline A}\Vol_{\OO\VV}$ and ${\bf 1}_{\OO\underline B}\Vol_{\OO\WW}$ respectively.
\end{lem}

These densities are denoted by $f^{(N)}_\VV(\be_0\,;\,\cdot)$ and $f^{(N)}_\WW(\be_0\,;\,\cdot)$ respectively\footnote{Recall $\eta$, and so $X_\VV$ and $X_\WW$, depend on $N$.}. We are going to prove that we have 
\begin{equation}
\label{AlmostIdentity}
q\bigl(\varpi_\VV(\be),\bfg^0\bigr)f^{(N)}_\WW(\be_0\,;\,\be) = f^{(N)}_\VV(\be_0\,;\,\be);
\end{equation}
we shall then get identity \eqref{FundRel} using the following lemma.

\begin{lem}
\label{Limites}
\begin{itemize}
   \item $f^{(N)}_\VV(\be_0\,;\,\be)\underset{N+\infty}{\longrightarrow} f_\VV(\be_0\,;\,\be)$.
   \item $f^{(N)}_\WW(\be_0\,;\,\be)\underset{N+\infty}{\longrightarrow} f(\be_0\,;\,\be)$.
\end{itemize}
\end{lem}

\smallskip

To proceed further and establish identity \eqref{AlmostIdentity}, we need to give some definitions. If $\VV$ and $\eta $ are chosen small enough, there exists a bundle isomorphism trivializing $\OO\mcU$:
$$
\psi : \RR^{1,d}\times SO_0(1,d)\rightarrow \OO\mcU.
$$
We shall denote by $(\zeta,g)$ the point $\psi^{-1}(\be)=\psi^{-1}\bigl((m,\bfg)\bigr)$; the set $Z$ will be the intersection of a small Euclidean ball of $\mathcal{M}_{1+d}(\RR)$, of center $g$, with $SO_0{(1,d)}$. Set $\mcA_Z := \OO A \cap\psi\bigl(\RR^{1,d}\times Z\bigr)$. The set $B_Z\subset\WW$ is the defined as the set of points of $\WW$ of the form $\exp_{\widehat m}(s\widehat T^0)$, with $(\widehat m,\widehat T)\in\mcA_Z$. The collection of all the $(m',\bfg')\in\OO B_Z$ where $m' = \exp_{\widehat m}(s\widehat T^0)$ and $\bfg'$ is the image of $\widehat T$ by parallell transport along the geodesic $\exp_{m'}(\cdot\,T^0)$ is denoted by $\mcB_Z$. Similarly, $A'_Z$ is defined as the set of points of $\VV$ of the form $\exp_{\widehat m}(s\widehat T^0)$, with $(\widehat m,\widehat T)\in\mcB_Z$. The collection of all the $(m',\bfg')\in\OO\VV$, where $m' = \exp_{\widehat m}(s\widehat T^0)$, $(\widehat m, \widehat T)\in\mcB_Z$, and $\bfg'$ is the image of $\widehat T$ by parallell transport along the geodesic $\exp_{\widehat m}(\cdot\,\widehat T^0)$, is denoted by $\mcA'_Z$. The range of $s$ in these definitions is restricted in such a way that the geodesics $\exp_\cdot\bigl(\cdot\,T^0\bigr)$ remain in $\mcU$. 

\smallskip

To prove identity \eqref{AlmostIdentity} we start from the inclusions 
$$
\bigl\{X_\VV\in\mcA_Z\bigr\}\subset\bigl\{X_\WW\in\mcB_Z\bigr\}\subset\bigl\{X_\VV\in\mcA'_Z\bigr\}
$$
to get the inequalities
$$
1\leq \frac{\PP_{\be_0}\bigl(X_\WW\in\mcB_Z\bigr)}{\PP_{\be_0}\bigl(X_\VV\in\mcA_Z\bigr)} \leq \frac{\PP_{\be_0}\bigl(X_\VV\in\mcA'_Z\bigr)}{\PP_{\be_0}\bigl(X_\VV\in\mcA_Z\bigr)},
$$
\textit{i.e.}
\begin{equation}
\label{IneqUn}
1\leq \frac{\int f_\WW^{(N)}(\be_0\,;\,\widehat\be){\bf 1}_{\mcB_Z}(\widehat\be)\Vol_{\OO\WW}(d\widehat\be)}{\int f_\VV^{(N)}(\be_0\,;\,\widehat\be){\bf 1}_{\mcA_Z}(\widehat\be)\Vol_{\OO\VV}(d\widehat\be)} \leq \frac{\PP_{\be_0}\bigl(X_\VV\in\mcA'_Z\bigr)}{\PP_{\be_0}\bigl(X_\VV\in\mcA_Z\bigr)}.
\end{equation}
We are going to obtain identity \eqref{AlmostIdentity} taking successively the supremum limit in the above inequalities, first as $A$ decreases to $\{m\}$, and then as $Z$ descreases to $\{g\}$. This final step rests on the following fact.

\begin{lem}
\label{LemmeLimite}
$\displaystyle{\underset{Z \searrow \{g\}}{\overline\lim}\;\underset{A \searrow \{m\}}{\overline\lim} \frac{\PP_{\be_0}\bigl(X_\VV\in\mcA'_Z\bigr)}{\PP_{\be_0}\bigl(X_\VV\in\mcA_Z\bigr)} = 1}$.
\end{lem}

\smallskip

This comes from the fact that the ratio of the $\Vol_{\OO\VV_{-\eta}}$-volume of the sets $\bigl\{\widetilde\be\in\OO\VV_{-\eta}'\,;\,X_\VV(\widetilde\be)\in\mcA_Z\bigr\}$ and $\bigl\{\widetilde\be\in\OO\VV_{-\eta}'\,;\,X_\VV(\widetilde\be)\in\mcA'_Z\bigr\}$ converges to $1$ as $A \searrow \{m\}$ and $Z\searrow \{g\}$. Recall $\sigma_\VV$ and $\sigma_\WW$ are the volume element induced by $q$ on $\VV$ and $\WW$ respectively. It remains to evaluate the ratio of the integals in equation \eqref{IneqUn} to get the conclusion; to that end we use the following fact.

\begin{lem}
\label{DernierLemme}
\begin{itemize}
   \item The limit $\displaystyle{\quad\underset{A\searrow \{m\}}{\lim}\,\frac{\sigma_\WW(B_Z)}{\sigma_\VV(A)}}$ exists, and
   \item there exists a positive function $c(Z)$ of $Z$, decreasing to $0$ as $Z$ decreases to $\{g\}$, and \st we have
$$
\bigl(1-c(Z)\bigr)q\bigl(\varpi_\VV(\be),\bfg^0\bigr) \leq \underset{A\searrow \{m\}}{\lim}\,\frac{\sigma_\WW(B_Z)}{\sigma_\VV(A)} \leq \bigl(1+c(Z)\bigr)q\bigl(\varpi_\VV(\be),\bfg^0\bigr)
$$
for $Z$ small enough.
\end{itemize}
\end{lem}
$\circ$ To see where this result comes from, write $\exp_m^\VV : T_m\VV \rightarrow \VV$ for the exponential map in $\VV$ at point $m$, and $\exp_m^\WW : T_m\WW \rightarrow \WW$ for the exponential map in $\WW$ at point $m$. We measure volumes in $T_m\VV$ and $T_m\WW$ using the (constant) volume elements $\Vol_m^\VV$ and $\Vol_m^\WW$ induced by $q$ on $T_m\VV$ and $T_m\WW$ respectively. Writing $A = \exp_m^\VV\bigl(\widetilde A\bigr)$ we have
\begin{equation}
\label{EstimateSurface1}
\frac{\sigma_\VV(A)}{\Vol_m^\VV\bigl(\widetilde A\bigr)}\underset{A\searrow \{m\}}{\longrightarrow} 1
\end{equation}
Associate to $\widetilde\bfg\in\OO_m\MMM\cap\psi\bigl(\RR^{1,d}\times Z\bigr)$ the set $\widetilde B_{\widetilde\bfg}\subset T_m\MMM$, image of $\widetilde A\subset T_m\VV$ by the projection map $T_m\MMM\rightarrow T_m\WW$ parallelly to $\widetilde\bfg^0$. We have on the one hand
\begin{equation*}
\frac{\Vol_m^\WW\left(\bigcup_{\widetilde\bfg}{\widetilde B}_{\widetilde\bfg}\right)}{\sigma_\WW(B_Z)} \underset{A\searrow \{m\}}{\longrightarrow} 1,
\end{equation*}
where the union is taken over all $\widetilde\bfg\in\OO_m\MMM\cap\psi\bigl(\RR^{1,d}\times Z\bigr)$, and on the other hand
\begin{equation*}
\frac{\Vol_m^\WW\bigl(\widetilde B_{\widetilde\bfg}\bigr)}{\Vol_m^\VV(\widetilde A)} = q\bigl(\varpi_\VV(\be),\widetilde\bfg^0\bigr).
\end{equation*}
Together with limit \eqref{EstimateSurface1} these two estimates imply lemma \ref{DernierLemme}. $\circ$

\smallskip

Decomposing $\mcB_Z$ into the union of its fibers: $\mcB_Z=: \displaystyle{\bigcup_{\widehat m\in B_Z}\mcB_Z^{\widehat m}}$, we can write the integral $\displaystyle{\int f_\WW^{(N)}(\be_0\,;\,\widehat\be)\,{\bf 1}_{\mcB_Z}(\widehat\be)\Vol_{\OO\WW}(d\widehat\be)}$ as 
$$
\int_\WW\left(\int_{\OO_{\widehat m}\MMM}f_\WW^{(N)}\bigl(\be_0\,;\,(\widehat m, \widehat\bfg)\bigr){\bf 1}_{\mcB_Z^{\widehat m}}(\widehat\bfg)\Vol_{\widehat m}(d\widehat\bfg)\right){\bf 1}_{B_Z}(\widehat m)\Vol_{\WW}(d\widehat m).
$$
A similar decomposition can be written for $\displaystyle{\int f_\VV^{(N)}(\be_0\,;\,\widehat\be)\,{\bf 1}_{\mcA_Z}(\widehat\be)\Vol_{\OO\VV}(d\widehat\be)}$ using the decomposition $\mcA_Z=: \displaystyle{\bigcup_{\widehat m\in B_Z}\mcA_Z^{\widehat m}}$ of $\mcA_Z$ into fibers:
$$
\int_\VV\left(\int_{\OO_{\widehat m}\MMM}f_\VV^{(N)}\bigl(\be_0\,;\,(\widehat m, \widehat\bfg)\bigr){\bf 1}_{\mcA_Z^{\widehat m}}(\widehat\bfg)\Vol_{\widehat m}(d\widehat\bfg)\right){\bf 1}_{A_Z}(\widehat m)\Vol_{\VV}(d\widehat m).
$$
Note that $\mcA_Z$ and $\mcB_Z$ have the same fiber at point $m$, namely $\mcB_Z^{m} = \OO_m\MMM\cap\psi\bigl(\RR^{1,d}\times Z\bigr)$. We get as a consequence of lemma \ref{DernierLemme} the following two inequalities:
\end{DemOuv}
\begin{equation*}
\bigl(1-c(Z)\bigr) q\bigl(\varpi_\VV(\be),\bfg^0\bigr)\frac{\int_{\OO_m\MMM} f_\WW^{(N)}\bigl(\be_0\,;\,(m,\bfg)\bigr){\bf 1}_{\mcB_Z^m}(\bfg)\Vol_{\OO_m\MMM}(d\bfg)}{\int_{\OO_m\MMM} f_\VV^{(N)}\bigl(\be_0\,;\,(m,\bfg)\bigr){\bf 1}_{\mcB_Z^m}(\bfg)\Vol_{\OO_m\MMM}(d\bfg)} \leq \underset{A \searrow \{m\}}{\underline\lim} \frac{\int f_\WW^{(N)}(\be_0\,;\,\widehat\be){\bf 1}_{\mcB_Z}(\widehat\be)\Vol_{\OO\WW}(d\widehat\be)}{\int f_\VV^{(N)}(\be_0\,;\,\widehat\be){\bf 1}_{\mcA_Z}(\widehat\be)\Vol_{\OO\VV}(d\widehat\be)}
\end{equation*}
$\quad\;\;\,$ and
\begin{equation*}
\underset{A \searrow \{m\}}{\overline\lim} \frac{\int f_\WW^{(N)}(\be_0\,;\,\widehat\be){\bf 1}_{\mcB_Z}(\widehat\be)\Vol_{\OO\WW}(d\widehat\be)}{\int f_\VV^{(N)}(\be_0\,;\,\widehat\be){\bf 1}_{\mcA_Z}(\widehat\be)\Vol_{\OO\VV}(d\widehat\be)} \leq \bigl(1+c(Z)\bigr) q\bigl(\varpi_\VV(\be),\bfg^0\bigr)\frac{\int_{\OO_m\MMM} f_\WW^{(N)}\bigl(\be_0\,;\,(m,\bfg)\bigr){\bf 1}_{\mcB_Z^m}(\bfg)\Vol_{\OO_m\MMM}(d\bfg)}{\int_{\OO_m\MMM} f_\VV^{(N)}\bigl(\be_0\,;\,(m,\bfg)\bigr){\bf 1}_{\mcB_Z^m}(\bfg)\Vol_{\OO_m\MMM}(d\bfg)}.
\end{equation*}
\begin{DemFerm}
Taking the supremum limit as $Z$ decreases to $\{g\}$ and using lemma \ref{DernierLemme} we obtain
$$
\underset{Z \searrow \{g\}}{\overline\lim}\,\underset{A \searrow \{m\}}{\underline{\overline\lim}} \frac{\int f_\WW^{(N)}(\be_0\,;\,\widehat\be){\bf 1}_{\mcB_Z}(\widehat\be)\Vol_{\OO\WW}(d\widehat\be)}{\int f_\VV^{(N)}(\be_0\,;\,\widehat\be){\bf 1}_{\mcA_Z}(\widehat\be)\Vol_{\OO\VV}(d\widehat\be)} = q\bigl(\varpi_\VV(\be),\bfg^0\bigr) \frac{f^{(N)}_\WW(\be_0\,;\,\be)}{f^{(N)}_\VV(\be_0\,;\,\be)}.
$$
As equation \eqref{IneqUn} together with lemma \ref{LemmeLimite} tells us that this supremum limit is equal to $1$, we conclude that
$$
q\bigl(\varpi_\VV(\be),\bfg^0\bigr)f^{(N)}_\WW(\be_0\,;\,\be) = f^{(N)}_\VV(\be_0\,;\,\be).
$$
Identity \eqref{FundRel} follows from lemma \ref{Limites} sending $N$ to infinity.
\end{DemFerm}

\medskip

The property of the one-particle distribution function emphasized in proposition \ref{InterpretationOnePart} will be used to prove the following fundamental theorem.

\begin{thm}
\label{LStarOnePartFunction}
We have $\displaystyle{L^*f(\be_0\,;\,\cdot) = 0}$ in $\OO\MMM\backslash\{\be_0\}$.
\end{thm}

The approach to relativistic Ornstein-Uhlenbeck process and (general) relativistic diffusions developped so far in the work of F. Debbasch and his co-authors relies entirely on a similar (manifestly covariant) transport equation, which is given as the fundamental object in their approach. This theorem provides a dynamical justification of this approach.

\begin{Dem}
As theorem \ref{LStarOnePartFunction} is of a local nature, we are going to take for each point $\be\neq\be_0$ a $\VV\in\mcalV_\be$ and work in the neighbourhood $\OO\mcU$ of $\be$ constructed in the proof of proposition/definition \ref{DefnOnePartDistr}  using the normal variation of $\VV$. We shall use here the same notations as there; notice in addition that in an expresion like $q\bigl(\varpi(\be),\bfg^0\bigr)$, the vector $\varpi(\be)$ will be seen as an element of $T_m\MMM$ rather than its horizontal lifting.

The beginning of the proof is exactly the same as in point \textit{1, \textbf{a)}} of the proof of proposition/definition \ref{DefnOnePartDistr}. We repeat it here to ease the reading ; it is quoted between the two stars $(*)$.

\medskip

\textit{1)} $(*)$ We shall suppose first that $\be_0$ belongs to $\OO\mcU$. If $\be_0=(m_0,\bfg^0)$ does not belong to the (closure of the) timelike past of $\mcU$, then no timelike path started from $m_0$ can ever hit $\VV$, so the function $f(\be_0\,;\,\cdot)$ is null in a neighbourhood of $\OO\VV$. As we are interested in what happens near $\OO\VV$, we shall make the hypothesis that $\be_0$ belongs to the timelike past of $\VV$. We shall suppose, \wlg, that $\be_0\in\VV_{-\eta}$. It will be fixed throughout this paragraph.

\smallskip

As the hitting times $H_\varep$ will be $\PP_{\be_0}$-almost surely finite under the preceding hypothesis, we can consider the re-parametrized process $\{\be_{H_\varep}\}_{\varep\in(-\eta,\eta)}$; it has generator $\overline\gamma^{-1}\,L$. We shall decompose this oeprator under the form
\begin{equation}
\label{DecompositionGenerator}
\forall\,\be=\phi_\varep(\widehat\be),\quad \frac{Lf}{\overline\gamma}(\be) = (\varpi f)(\be) + \widehat L(f\circ \phi_\varep)\,(\widehat\be) = (\varpi f)(\be) + \bigl(\overline L f\bigr)(\be),
\end{equation}
where $\widehat L$ is a second order differential operator on $\OO\VV$, and where, as a consequence, $\overline L$ acts only on $\OO\VV_\varep$. Now, define the $\OO\VV$-valued process $\bigl\{\widehat\be_\varep\bigr\}_{\varep\in(-\eta,\eta)} := \bigl\{\phi_\varep^{-1}(\be_{H_\varep})\bigr\}_{\varep\in(-\eta,\eta)}$ and denote by $\widehat\ell_\varep$ its time-dependent generator.(*) This operator is seen to be hypoelliptic, so the \rv $\widehat\be_\varep$ has for any $\varep\in (-\eta,\eta)$ a smooth density $\widehat\rho(\be_0\,;\,\varep,\cdot)$ \wrt $\Vol_{\OO\VV}$ which satisfies the equation 
$$
\forall\varep\in(-\eta,\eta),\quad\quad \partial_\varep \widehat\rho(\be_0\,;\,\varep,\cdot) = \widehat\ell_\varep^{*\OO\VV}\widehat\rho(\be_0\,;\,\varep,\cdot).
$$ 
$\widehat\ell_\varep^{*\OO\VV}$ stands here for the $\LL^2(\Vol_{\OO\VV})$-dual of $\widehat\ell_\varep$. Let us now denote by $\Vol_{\OO\VV}^\varep$ the pull-back on $\OO\VV$ by $\phi_\varep$ of the measure $\Vol_{\OO\VV_\varep}$ on $\OO\VV_\varep$, and denote by $G_\varep$ its density \wrt $\Vol_{\OO\VV}$. Then $\widehat\be_\varep$ has a density $\displaystyle{\widehat\mu_\varep(\be_0\,;\,\varep,\cdot) = \frac{\widehat\rho_\varep(\be_0\,;\,\varep,\cdot)}{G_\varep}}$ \wrt $\Vol_{\OO\VV}^\varep$; it satisfies the equation
\begin{equation}
\label{HeatCharts}
\partial_\varep \widehat\mu_\varep(\be_0\,;\,\varep,\cdot) + \frac{\partial_\varep G_\varep}{G_\varep}\widehat\mu_\varep(\be_0\,;\,\varep,\cdot) = \widehat\ell_\varep^{*\OO\VV;\varep}\widehat\mu_\varep(\be_0\,;\,\varep,\cdot).
\end{equation}
We have here \,$\widehat\ell_\varep^{*\OO\VV;\varep} g = \frac{\widehat\ell_\varep^{*\OO\VV}(G_\varep g)}{G_\varep}\,$ for any smooth function $g$. Denote by $\mu(\be_0\,;\,\varep,\cdot)$ the density of $\be_{S_\varep}$ \wrt $\Vol_{\OO\VV_\varep}$, and consider $\mu$ and $G$ as functions of $\varep$ and $\be\in\VV_\varep$, \textit{i.e.} consider them as functions defined on the open set $\mcU$. Then, equation \eqref{HeatCharts} can be written
\begin{equation}
\label{HeatCharts2}
\varpi\mu(\be_0\,;\,\cdot) + \frac{\varpi\,G}{G}\mu(\be_0\,;\,\cdot) = \overline{L}^{*\OO\VV_\varep}\mu(\be_0\,;\,\cdot).
\end{equation}
The operator $\overline L$ has been introduced in equation \eqref{DecompositionGenerator}. It is useful at that stage to remark that we have\footnote{Recall that $\Vol$ is the Liouville measure on $\OO\MMM$ and that the $^*$-operation is the $\LL^2(\Vol)$-dual operation.} 
$$
\overline{L}^{*\OO\VV_\varep} = \overline{L}^* 
$$
as a consequence of the change of variable formula, and since we have a \textit{normal} variation of $\VV$. The following lemma is needed to make the final step.

\begin{lem}
We have for any smooth function $f$
$$
\varpi^* f + \varpi f +\frac{\varpi\,G}{G} f = 0.
$$
\end{lem}

\begin{SousDem}
As above, this is consequence of the change of variable formula and the fact that we have a normal variation of $\VV$. We have, for any smooth function $\varphi$ with compact support,
\begin{equation*}
\begin{split}
\int\left(\varpi^* f\right)(\be)\,&\varphi(\be)\Vol(d\be) = \int f(\be)\,(\varpi\,\varphi)(\be)\Vol(d\be) = \int f(\varep,\widehat\be)\,(\partial_\varep\varphi)(\varep,\widehat\be)\,G_\varep(\widehat\be)\,\sigma_0(d\widehat\be)\,d\varep \\ 
&= -\int (\partial_\varep f)(\varep,\widehat\be)\,\varphi(\varep,\widehat\be)\,G_\varep(\widehat\be)\,\sigma_0(d\widehat\be)\,d\varep - \int (f\varphi)(\varep,\widehat\be)\,\partial_\varep G_\varep(\widehat\be)\,\sigma_0(d\widehat\be)\,d\varep \\
&= -\int\left(\varpi f + \frac{\varpi\,G}{G}\right)(\be)\,\varphi(\be)\,\Vol(d\be).
\end{split}
\end{equation*}
\end{SousDem}

As a consequence of this lemma we can use the decomposition given in equation \eqref{DecompositionGenerator} to write equation \eqref{HeatCharts2} as
$$
L^*\left(\frac{\mu(\be_0\,;\,\cdot)}{\overline\gamma}\right) = 0.
$$
Proposition \ref{InterpretationOnePart} enables to conclude that $L^*f(\be_0\,;\,\cdot) = 0$ in $\mcU$.

\medskip


\textit{2)} To deal with the general case where $\be_0$ does not belong to $\OO\mcU$, denote by $\VV'\subset\VV$ an open subset of $\VV$ such that any timelike path hitting $\VV'$ hits $\VV_{-\eta}$ before. Such a manifold $\VV'$ will exist provided $\eta$ is small enough. For a small enough $\delta>0$, the set $\phi_{(-\delta,\delta)}(\VV') := \bigl\{\phi_\varep(\widehat m)\in\MMM\,;\,\widehat m\in\VV' \textrm{ and }|\varep|<\delta\bigr\}$ will have the property that any timelike path hitting it hits $\VV_{-\eta}$ before. Set 
$$
\mcU' := \OO\phi_{(-\delta,\delta)}(\VV');
$$
this is an open set of $\OO\MMM$. To prove theorem \ref{LStarOnePartFunction} on $\mcU'$, it suffices to remark that for $\be\in\mcU'$
, we have
$$
f(\be_0\,;\,\be) = \EE_{\be_0}\bigl[f(\be_H\,;\,\be){\bf 1}_{H<\infty}\bigr],
$$
where $H$ is the hitting time of $\VV_{-\eta}$. The previous part of the proof applies to each function $f(\be_H\,;\,\cdot)$. It follows then from the above identity that $f(\be_0\,;\,\cdot)_{\big| \mcU'}$ is an $L^*$-harmonic function, as a mean of $L^*$-harmonic functions.
\end{Dem}

\subsubsection{$L$-harmonic functions}
\label{SubSubSectionHarmonicFunctions}
We shall see in section \ref{BouquetFinal} an important application of theorem \ref{LStarOnePartFunction} in relation with statistical irreversibility. Before turning oursleves to that side, we would like to stress in this section the importance that theorem \ref{LStarOnePartFunction} might have from a geometrical point of view. To that end, we shall investigate its meaning in the study of the $(0,\be_.)$-diffusion of Franchi and Le Jan. As emphasized after definition \ref{DefnDiff}, this $(V,\ff)$-diffusion is the only process of this class determined entirely by the geometric background $(\MMM,q)$; this property gives it a special status. Its generator is 
$$
L = H_0 + \frac{1}{2}\sum_{i=1}^d V_i^2.
$$
We shall call a $\mcC^2$ function on $\OO\MMM$ satisfying the relation $Lf=0$ an \textit{\textbf{L-harmonic function}}. The class of bounded $L$-harmonic functions and the asymptotic behaviour of the $(0,\be_.)$-diffusion are two faces of the same object: the boundary at infinity of the manifold $(\MMM,q)$. 

\paragraph{a) Ideal boundaries of manifolds and invariant $\sigma$-algebra.} Let us illustrate this correspondence recalling what happens to \bm on some special Riemannian manifolds; as the $(0,\be_.)$-diffusion, \bm is entirely determined by the geometric environment. Suppose $(\MMM,q)$ is a simply connected Cartan-Hadamard manifold: it is a Riemannian manifold, diffeomorphic to some $\RR^n$, with curvature bounded by two negative constants. The exponential polar co-ordinates $(r,\theta)\in\RR_+\times\mathbb{S}^{d-1}$ associated with any point provide global co-ordinates on $\MMM$. These manifolds have the property any sequence of balls $\{B_i\}_{i\geq 0}$, with constant radius, whose centers leave any compact, appear uniformly small when seen from within a compact set: For any compact set $K$ and given any $\varep>0$, there exists an index $i_\varep$ \st for any $i\geq i_\varep$ and any point $m\in K$, whose accociated system of polar co-ordinates is denoted by $(r\,\theta)$, any point of $B_i$ has polar angle $\theta$ contained in a region of $\mathbb{S}^{d-1}$ of diameter no greater than $\varep$.  These manifolds also enjoy the following property: Given any geodesic $\{\gamma_t\}_{t\geq 0}$ and any point $m'$, there exists a unique geodesic $\{\gamma'_t\}_{t\geq 0}$ started from $m'$ \st the distance between $\gamma'_t$ and $\gamma_t$ remains bounded. These properties motivate the introduction of a compactification of $\MMM$, homeomorphic to $\mathbb{S}^{d-1}$, and where a path converges to a point of tboundary if its polar angle converges in any polar system of co-ordinates. Any geodesic converges to some point of the boundary\footnote{See for instance the article \cite{Anderson} of Anderson, or chapter $8$ of the book \cite{Pinsky} of R. Pinsky.}.

From a probabilistic point of view, we can investigate the far end of a manifold looking at what happens to \bm $\{w_t\}_{t\geq 0}$ on $(\MMM,q)$ as time goes to $\infty$. To that end, we define the \textit{\textbf{invariant $\sigma$-algebra}} of the process as being generated by the real-valued functionals $F(w_.)$ depending only on the asymptotic beahaviour of $w_.$, \textit{i.e.} satisfying the identity $F\bigl(\{w_t\}_{t\geq 0}\bigr) = F\bigl(\{w_{s+t}\}_{t\geq 0}\bigr)$ for any $s\geq 0$. 

Back in the above Cartan-Hadamard manifold $(\MMM,q)$, pick a point $m\in\MMM$ and denote by $(r_t,\theta_t)$ the $m$-polar co-ordinates of $w_t$. It can be proved\footnote{See for instance the pioneering article \cite{Anderson} of Anderson, or the article \cite{Kifer} of Y. Kifer.} that $\{w_t\}_{t\geq 0}$ converges $\PP_{w_0}$-\as to some random point of $\partial\MMM$, characterized by the fact that $\theta_t\rightarrow \theta_\infty\in\mathbb{S}^{d-1}$. It can also be shown that the invariant $\sigma$-algebra is $\PP_{w_0}$-indistiguishable from the algebra generated by $\theta_\infty$. This fact gives a probabilistic meaning to $\partial\MMM$, or, conversely, gives a geometric meaning to the invariant $\sigma$-algebra of Brownian motion.

The situation appears to be similar, though subtler, in the Lorentzian framework of Minkowski space. Recall the \textbf{\textit{causal boundary}} ${\bf C}$ of $\RR^{1,3}$ is the ideal boundary of $\RR^{1,3}$ characterized by the property that two timelike paths $\{\gamma_t\}_{t\geq 0}$ and $\{\gamma'_t\}_{t\geq 0}$ converge to the same boundary point iff they have the same chronological past: $I^-(\gamma) = I^-(\gamma')$. The following theorem has been proved in the articles \cite{BailleulPoisson} of I. Bailleul and \cite{BailleulRaugi} of I. Bailleul and A. Raugi. It holds for any starting point $\be_0$ of the $(0,\be_.)$-diffusion.

\begin{thm}[\cite{BailleulPoisson}, \cite{BailleulRaugi}]
 \begin{itemize}
   \item The $\RR^{1,d}$-part $\{m_s\}_{s\geq 0}$ of the $(0,\be_.)$-diffusion converges $\PP_{\be_0}$-almost surely to some random point $m_{\infty}$ of ${\bf C}$.
   \item The $\sigma$-algebra generated by $m_{\infty}$ coincides with the tail $\sigma$-algebra of $\{\xi_s\}_{s\geq 0}$, up to $\PP_{\be_0}$-null sets. 
\end{itemize}
\end{thm}

So, we can find back the causal boundary in the probabilistic invariant $\sigma$-algebra. This is a nice feature that might help clarify geometrically more complicated situations, giving a simple probabilistic picture of what happens. J. Franchi has for example undertaken in \cite{Franchi} the study of the $(0,\be_.)$-diffusion in G\"odel's spacetime. This space has a trivial causal boundary, reduced to one point. Yet, he has been able to prove that the invariant $\sigma$-algebra of the process is not trivial. This suggested, in return, the definition of a purely geometric boundary.

\paragraph{b) Poisson and Martin boundaries.} 
The link between geometry and probability illustrated above is complemented by the existing link between invariant $\sigma$-algebra on the one hand and the set of bounded $L$-harmonic functions on the other hand\footnote{Consult for instance chapter $8$ of the book \cite{Pinsky} of R. Pinsky for the Riemannian case, and proposition $8$ in the article \cite{BailleulPoisson} for the hypoelliptic situation appearing in the study of the $(0,\be_.)$-diffusion in Minkowski space.}. It is equivalent to determine one or the other. The set of bounded $L$-harmonic functions is called the \textbf{\textit{Poisson boundary of $(L,\MMM)$}}. So, the Poisson boundary of $L$, the invariant $\sigma$-algebra and the geometry at infinity of $(\MMM,q)$ may be seen as three faces of a same object.

\smallskip

Let us give a last picture of the Riemannian/Brownian situation. We shall get a clearer image looking at any elliptic smooth second order differential operator $L_0$ on a connected (relatively compact) open set $D$ of $\RR^n$. Recall the \textbf{\textit{Martin boundary of $(L_0,D)$}} is the collection of non-negative $L_0$-harmonic functions on $D$(\footnote{This set contains the Poisson boundary of $L_0$ in $D$.}). Martin gave in \cite{Martin} a methof to construct this set and proved that any non-negative $L_0$-harmonic function can be uniquely represented as the barycenter of a finite measure on the set of extreme points of his boundary. This construction is now well understood from a probabilistic point of view (see for instance chapter $7$ of Dynkin's book \cite{Dynkin}). Let us briefly describe it. Denote by $G(x,y)$ the Green kernel of $L_0$ in $D$, and define the function 
$$
K(x,y) = \frac{G(x,y)}{G(x_0,y)}, \quad x\in D, y\in D\backslash\{x_0\},
$$
where $x_0$ is a fixed point. Observe that the function $K(\cdot,y)$ is $L_0$-harmonic on $D\backslash\{y\}$, for any $y\in D\backslash\{x_0\}$. It follows that we shall construct $L_0$-harmonic functions on $D$ sending $y$ to the boundary of $D$, provided the limit $\lim_y K(\cdot,y)$ exists. Martin's boundary is made up of all the functions obtained that way. A sequence $\{y_n\}_{n\geq 0}$ of points of $D$ leaving every compact $D$, and \st the function $K(\cdot,y_n)$ converges\footnote{Uniformly locally on compact subsets of $D$.}, is called a \textbf{\textit{fundamental sequence for $(L_0,D)$}}. In short, Martin's theory asserts that the knowledge of fundamental sequences is equivalent to the knowledge of the set of non-negative $L_0$-harmonic functions. 

Pinsky gave in \cite{Pinsky93} a probabilistic proof of a useful characterisation of fundamental sequences, known from potential theoretists before\footnote{See the notes of chapter $7$ in \cite{Pinsky}.}. $L_0^*$ will denote the $\LL^2(\LEB)$-dual of $L_0$ and $\{\widetilde{\PP}_x\}_{x\in D}$ the laws of the diffusion $\{X_t\}$ in $D$ with generator $L_0^*$. We shall denote by $\zeta$ its exit time from $D$, and, given a compact subset $U$ of $D$, we shall denote by $H_U$ the hitting time of $U$ by $\{X_t\}_{0\leq t <\zeta}$.

\begin{thm}[Pinsky\cite{Pinsky93}]
The sequence $\{y_n\}_{n\geq 0}$ is fundamental for $(L_0,D)$ iff, for any smooth compact subset $U$ of $D$, the sequence of conditional distributions $\bigl\{\widetilde{\PP}_{y_n}\bigl(X_{H_U}\in \cdot\,\big|\,H_U<\zeta \bigr)\bigr\}_{n\geq 0}$ converges.
\end{thm}

\paragraph{c) A conjecture.} 
Theorem \ref{LStarOnePartFunction} bringing into play $L^*$ and hitting distributions, through $f(\be_0\,;\,\cdot)$, it is now time to examine it. Note that 
$$
L^* = -H_0 + \frac{1}{2}\sum_{i=1}^dV_i^2.
$$
So, $L^*$ is the genrator of an $\OO\MMM$-valued diffusion analogue to the $(0,\be_.)$-process, except that the speed $\frac{dm_s}{ds} = -\bfg^0_s$ is past-directed. Call it \textbf{\textit{$(0,\bel_.)$-diffusion}} and denote by $\overset{\leftarrow}{\PP}_{\be_0}$ its law when started from $\be_0\in\OO\MMM$. It is clear from its construction that the paths of $\bel_.$ started from $\be_0$ take values in the chronological past $I^-(\be_0)$ of $\be_0$. Any open set of $I^-(\be_0)$ is visited by the process with positive probability, and any spacelike hypersurface $\VV$ is hit with positive probability. As noted in proposition \ref{InterpretationOnePart}, these hitting distributions are determined by the one-particle distribution function of the $(0,\bel_.)$-process. We shall denote it by $\overset{\leftarrow}{f}(\be_0\,;\,\cdot),\,\be_0\in\OO\MMM$. Theorem \ref{LStarOnePartFunction} can be restated as follows.

\begin{figure}[h!]
\begin{center}
\input{ReversedDiffusion.pstex_t} 
\end{center}
\caption{$(0,\be_.)$-diffusion and $(\bel_.,0)$-diffusion.}
\end{figure}

\begin{thm}
The function $\overset{\leftarrow}{f}(\be_0\,;\,\cdot)$ is $L$-harmonic on $\OO\MMM\backslash\{\be_0\}$.
\end{thm}
It is tempting, after reading paragraph \textbf{b)}, to renormalize $\overset{\leftarrow}{f}(\be_0\,;\,\cdot)$ and try to get possibly non-null $L$-harmonic functions sending the singularity $\be_0$ to infinity. This could be made looking at 
$$
\frac{\overset{\leftarrow}{f}(\be_0\,;\,\cdot)}{\overset{\leftarrow}{f}(\be_0\,;\,\bfc)}
$$
for some $c\in\OO\MMM$, or 
$$
(m,\bfg)\in\OO\MMM\mapsto\frac{\overset{\leftarrow}{f}\bigl(\be_0\,;\,(m,\bfg)\bigr)}{\int_{\OO_m\MMM} \overset{\leftarrow}{f}\bigl(\be_0\,;\,(m,\bfg')\bigr)\Vol_m(d\bfg')}.
$$
As emphasized in proposition \ref{InterpretationOnePart}, the use of the second ratio essentially amounts to look at the convergence of the hitting distributions $\PP_{\be_0}\bigl(\overset{\leftarrow}{\be}_{H_\VV}\in\cdot\,\big|\,H_\VV<\infty\bigr)$ of $\overset{\leftarrow}{\be}_.$ on any spacelike (smooth) hypersurface $\VV$. All this brings us to conjecture the following equivalence.

\begin{conj}
The following statements are equivalent.
\begin{enumerate}
   \item The sequence $\{\be_n\}_{n\geq 0}$ is fundamental for $L$ in $(\MMM,q)$.
   \item For any spacelike smooth hypersurface $\VV$ of $\MMM$ the sequence of conditional distributions 
$$
\overset{\leftarrow}{\PP}_{\be_n}\bigl(\overset{\leftarrow}{\be}_{H_\VV}\in\cdot\,\big|\,H_\VV<\infty\bigr)
$$ 
converges. 
   \item For any $\bfc\in\OO\MMM$, the sequence of $L$-harmonic functions\footnote{Given a compact set $K$ and a sequence $\{\be_n\}_{n\geq 0}$ of points of $\OO\MMM$ leaving every compact, the function $\overset{\leftarrow}{f}(\be_n\,;\,\cdot)$ is well defined on $K$ for $n$ large enough.} $\left\{\frac{\overset{\leftarrow}{f}(\be_n\,;\,\cdot)}{\overset{\leftarrow}{f}(\be_n\,;\,\bfc)}\right\}_{n\geq 0}$ converges uniformly on compact subsets of $I^-(\bfc)$.
\end{enumerate}
\end{conj}

\noindent This fact would explain why the causal boundary of $(\MMM,q)$ is likely to appear in the picture. In order for the conditional distributions $\overset{\leftarrow}{\PP}_{\be_n}\bigl(\overset{\leftarrow}{\be}_{H_\VV}\in\cdot\,\big|\,H_\VV<\infty\bigr)$ to converge, the support of each of these probabilities has to converge, for any spacelike hypersurface $\VV$. This cannot happen unless the chronological past $I^-(\be_n)=I^-\bigl((m_n,\bfg_n)\bigr)$ of $\be_n$ converges, \textit{i.e.} unless the sequence $\{m_n\}_{n\geq 0}$ has a limit in the causal boundary of $(\MMM,q)$. Yet, the study of G\"odel's spacetime by J. Franchi in \cite{Franchi} has made it clear that this geometric boundary might not be appropriate to describe the Poisson or the Martin boundary in some situations. Note, yet, that the above analysis using the one-particle distribution function does not apply in this non-strongly causal spacetime; no good definition of one-particle distribution function is available at the moment in such a framework. 

Last, we should oppose the difficulty of this problem on the \textit{large scale structure} of $(\MMM,q)$ to the previously mentionned fact that the \textit{local geometry of spacetime} can be recovered looking at the pathwise behaviour of the $(0,\be_.)$-process. Complications come from infinity... We shall come back to the above conjecture in a near future.

\subsection{H-theorem}
\label{BouquetFinal}

We give in this last section a proof of the analogue of the H-theorem first proved in \cite{ROUP02} for the R.O.U.P. in Minkowski space, as defined there through a 'Kolmogorov equation'. It has been then extended in \cite{ROUPHThm} to the R.O.U.P. (as defined in \cite{ROUPCurved}) in any Lorentzian manifold, and finally in \cite{ROUPUnifying1} to a larger class of 'diffusions' in Minkowski space. We deal here with the general case of $(V,\ff)$-diffusions in any Lorentzian manifold. We refer to the articles \cite{ROUP02}, \cite{ROUPHThm} and \cite{IsraelBoltzmann} for physical motivations.

\medskip

Let $U\subset\OO\MMM$ be a relatively compact open set and $f$ and $g$ be two positive smooth functions on $U$ satisfying the relations $L^*f=L^*g=0$(\footnote{These functions could for example be of the form $f(\be_0\,;\,\cdot)$ and $f(\be_0'\,;\,\cdot)$ if the strong causality assumption on $(\MMM,q)$ is satisfied.}). We shall denote by $Y$ a continuous unit vector field on $U$; define the function $\rho : \OO U \rightarrow (0,+\infty)$, $(m,\bfg)\mapsto q(Y_m,\bfg^0)$. 

\smallskip

We shall make the following assumptions on $f$ and $g$; they are sufficient to ensure the existence of the integrals below, and to differentiate them.
\begin{itemize}
   \item $\ln \frac{f}{g}$ is bounded.
   \item There exists positive constants $C$ and $\epsilon$ \st $f$ and its first and second derivatives are uniformly bounded by $Ce^{-\rho^{1+\epsilon}(m,\bfg)}$ in $U$.
\end{itemize}

Define now on $U$ the vector field 
$$
X(m) = -\int_{\OO_m\MMM}\bfg^0\,f(m,\bfg)\ln\frac{f(m,\bfg)}{g(m,\bfg)}\,\Vol_m(d\bfg).
$$
The main result of this section is the following theorem; no assumption on the geometry of space or on the data $V, \ff$ is needed.

\begin{thm}[H-theorem]
\label{ThmDivergence}
We have $\emph{div} X \geq 0$ for any two $L^*$-harmonic functions $f$ and $g$, and $X$ defined as above.
\end{thm}

We shall begin the proof of this theorem proving the following lemma\footnote{Compare with the Appendix to the article \cite{ROUPHThm} of F. Debbasch and M. Rigotti.}.

\begin{lem}
\label{LemUtile}
Given any (good) smooth function $h$ on $\OO\MMM$ set
$$
X(m) = \int_{\OO_m\MMM}\bfg^0\,h(m,\bfg)\,\Vol_m(d\bfg).
$$
Then, 
$$
(\emph{div}X)(m) = \int_{\OO_m\MMM}(H_0 h)(m,\bfg)\,\Vol_m(d\bfg).
$$
\end{lem}

\begin{Dem}
Given a $\mcC^1$ path $\gamma$ in $\MMM$ and two time $s,t$ we shall denote by $T^{\gamma}_{s\leftarrow t} : T_{\gamma_t}\MMM \rightarrow T_{\gamma_s}\MMM$ the parallel transport operation along the path $\{\gamma_r\}_{r\in [s,t]}$. It is an isometry between the two tangent spaces. We shall denote by $\nabla$ the Levi-Civita connection on $(\MMM,q)$. Recall that we have
$$
\nabla_{\dot\gamma_0} X = \underset{s\rightarrow 0}{\lim} \frac{T^\gamma_{0\leftarrow s}X_{\gamma_s} - X_{\gamma_0}}{s}.
$$
Recall also that the divergence of $X$ is the (Lorentzian) trace of the map $\nabla_. X$. It means that given any choice of orthonormal frame $\overline\bfg$ of $T_m\MMM$, the sum 
$$
\sum_{i=0}^d q(\overline\bfg^i,\overline\bfg^i)q(\nabla_{\overline\bfg^i}X,\overline\bfg^i)
$$
is independent of $\overline\bfg\in\OO_m\MMM$; this is, by definition, $\bigl(\textrm{div}X\bigr)(m)$. Last, recall that the vector field $H_0$ is defined as the generator of the lift to $\OO\MMM$ of the geodesic flow on $\HH\MMM$. Its dynamics $\bigl\{(m_s,\bfg_s)\bigr\}$ is determined by the condition $\frac{dm_s}{ds}=\bfg^0_s$ and the fact that $\bfg_s$ is parallelly transported along the path $\{m_s\}$.

Choose now a frame $\overline\bfg\in\OO_m\MMM$ and a path $\gamma^i$ \st $\gamma^i(0)=m$ and $\dot\gamma^i(0) = \overline\bfg^i$. Since parallel transport is an isometry, we can write
$$
X_{\gamma^i_s} = \int_{\OO_{\gamma_s}\MMM} \widehat\bfg^0\,h(\gamma^i_s,\widehat\bfg)\,\Vol_{\gamma^i_s}(d\widehat\bfg) = \int_{\OO_m\MMM} T^{\gamma^i}_{s\leftarrow 0}\bfg^0\,h\bigl(\gamma^i_s,T^{\gamma^i}_{s\leftarrow 0}\bfg\bigr)\,\Vol_m(d\bfg),
$$
so we have
$$
T^{\gamma^i}_{0\leftarrow s}X_{\gamma^i_s} = \int_{\OO_m\MMM} \bfg^0\,h\bigl(\gamma^i_s,T^{\gamma^i}_{s\leftarrow 0}\bfg\bigr)\,\Vol_m(d\bfg)
$$
and 
$$
\frac{T^{\gamma^i}_{0\leftarrow s}X_{\gamma^i_s} - X_m}{s} =  \int_{\OO_m\MMM} \bfg^0\;\frac{h\bigl(\gamma^i_s,T^{\gamma^i}_{s\leftarrow 0}\bfg\bigr)-h(m,\bfg)}{s}\,\Vol_m(d\bfg).
$$
Send $s$ to $0$ and sum over $i$ to get the result:
$$
(\textrm{div}X)(m) = \sum_{i=0}^d q(\overline\bfg^i,\overline\bfg^i)q(\nabla_{\overline\bfg^i}X,\overline\bfg^i)  = \int_{\OO_m\MMM}(H_0 h)(m,\bfg)\,\Vol_m(d\bfg).
$$
\end{Dem}

With this lemma in hand we can prove theorem \ref{ThmDivergence}. Recall that $L = H_0 + V + \frac{1}{2}V_iB^{ij}V_j$ and $L^* = -H_0 + V^* + \frac{1}{2}V_iB^{ij}V_j$.

\begin{DemOuv}
First, use lemma \ref{LemUtile} to write
$$
-\textrm{div}X = \int H_0\left(f\ln\frac{f}{g}\right)\,\Vol(d\bfg) = \int (H_0 f)\,\ln\frac{f}{g}\,\Vol(d\bfg) + \int \Bigl(H_0 f - \frac{f}{g}H_0 g\Bigr)\,\Vol(d\bfg).
$$
Use then the relations $L^*f=L^*g=0$ to get
\end{DemOuv}
\begin{equation*}
\begin{split}
&-\textrm{div} X = \int\textstyle{\Bigl(V^* f +\frac{1}{2}V_i\bigl(B^{ij}V_j(\la f)\bigr)\Bigr)\Bigl(\ln \frac{f}{g}+1\Bigr)\,\Vol(d\bfg)} - \int \textstyle{\frac{f}{g}\Bigl(V^* g +\frac{1}{2}V_i\bigl(B^{ij}V_j(\la g)\bigr)\Bigr)\,\Vol(d\bfg)} \\
&= \int\textstyle{\left((V^*f)\Bigl(\ln \frac{f}{g}+1\Bigr) - \frac{f}{g}\,V^* g\right)\,\Vol(d\bfg)} +\frac{1}{2}\int\textstyle{\left(\Bigl(\ln \frac{f}{g}+1\Bigr) V_i\bigl(B^{ij}V_j(\la f)\bigr) - \frac{f}{g} V_i\bigl(B^{ij}V_j(\la g)\bigr)\right)\,\Vol(d\bfg)}.
\end{split}
\end{equation*}
\begin{DemFerm}
Integrating by parts and using the relation $V\bigl(\ln \frac{f}{g}\bigr) = \frac{g}{f}V\bigl(\frac{f}{g}\bigr)$, the first integral is seen to be equal to 
$$
\int \left(g V\Bigl(\frac{f}{g}\Bigr) - V\Bigl(\ln\frac{f}{g}+1\Bigr)\right)\,\Vol(d\bfg) = 0.
$$
Recall $V_i^* = -V_i$. Use integration by parts in the second integral and the relation $V_j(\la f) = V_j(\la g)\frac{f}{g} + \la g V_j\Bigl(\frac{f}{g}\Bigr)$, to get
\begin{equation}
\begin{split}
-\textrm{div} X &= \frac{-1}{2} \int V_i\Bigl(\ln \frac{f}{g}\Bigr)B^{ij}V_j(\la f)\,\Vol(d\bfg) + \frac{1}{2} \int V_i\Bigl(\frac{f}{g}\Bigr)B^{ij}V_j(\la g)\,\Vol(d\bfg) \\
&= \frac{-1}{2} \int \left(\frac{g}{f}V_i\Bigl(\frac{f}{g}\Bigr)B^{ij}V_j(\la f)  - V_i\Bigl(\frac{f}{g}\Bigr)B^{ij}V_j(\la g)\right)\,\Vol(d\bfg) \\
&= \frac{-1}{2} \int \left\{V_i\Bigl(\frac{f}{g}\Bigr)B^{ij}\left(V_j(\la g)+\frac{\la g^2}{f} V_j\Bigl(\frac{f}{g}\Bigr)\right) - V_i\Bigl(\frac{f}{g}\Bigr)B^{ij}V_j(\la g)\right\}\,\Vol(d\bfg) \\
&= \frac{-1}{2} \int \frac{\la g^2}{f} V_i\Bigl(\frac{f}{g}\Bigr) B^{ij}V_j\Bigl(\frac{f}{g}\Bigr)\,\Vol(d\bfg).
\end{split}
\end{equation}
We get the conlustion from the non-negativeness of the matrix $B = \bigl(A^{-1}\bigr)^*A^{-1}$.
\end{DemFerm}

$\bullet$ Attention should be paid to the range of application of theorem \ref{ThmDivergence}. It seems tempting, indeed, in a strongly causal spacetime, to apply it to functions of the form $f(\be_0\,;\,\cdot)$ and $f(\be'_0\,;\,\cdot)$. However, today's state of art is far from being sufficient to provide estimates on these functions good enough to ensure that hypothese like those made at the beginning of the section hold\footnote{The boundedness hypothesis on $\frac{f}{g}$ is even most likely to be untrue for such functions.}. This is a difficult topic where the non-ellipticity of the operator $L$ complicates everything. As a first step towards such results, it would be interesting to determine small time estimates of its heat kernel; known results are unsufficient to answer this question.

\medskip

$\bullet$ Theorem \ref{ThmDivergence} proves that the flow of the vector field $X$ is volume increasing. It is not clear how one should interpret this result from a physical point of view when $\MMM$ is different from $\RR^{1,d}$. In this special case, choose a rest frame $\bfg\in SO(1,d)$ and denote by $t$ its associated time. Then, the integral of $X$ over any hyperplane of constant time is an increasing function of $t$ (provided $X$ is equal to $0$ at space infinity). This fact justifies that we should call theorem \ref{ThmDivergence} an $H$-theorem in that case. Things are less clear in any Lorentzian manifold, where time does not exist globally.

\smallskip

Things are even less satisfying from an information theoretic point of view. Recall that the relative entropy of a probability $\PP$ \wrt another probability $\QQ$ is infinite if $\PP$ is not absolutely continuous \wrt $\QQ$, and equal to 
$$
H(\PP\,;\,\QQ) = \EE_{\PP}\left[\ln \frac{d\PP}{d\QQ}\right],
$$
if  $\PP$ is absolutely continuous \wrt $\QQ$. We write $\EE_{\PP}$ for the expectation operator associated with $\PP$. Relative entropy is always non-negative, as is clear from the inequality $a\ln\frac{a}{b} \geq a-b$. 

Suppose for clarity that $\PP$ and $\QQ$ are probabilities on $[0,1]$(\footnote{This is not a serious restriction as any probability on a Borel space is isomorphic to a probability measure on $[0,1]$. This class of spaces is large enough to encompass most of the useful situations. See the Appendix of the book \cite{DynkinYushkevich} of Dynkin and Yushkevich.}). Let $X_1,X_2,...$ be i.i.d. random variables, with common law $\PP$ or $\QQ$. Then, given any numbers $x_1,...,x_n$ in $[0,1]$, we have 
$$
"\frac{\PP^{\otimes n}(X_1=x_1,..., X_n=x_n)}{\QQ^{\otimes n}(X_1=x_1,..., X_n=x_n)}" := \frac{d\PP^{\otimes n}}{d\QQ^{\otimes n}}(x_1,...,x_n) = \frac{d\PP}{d\QQ}(x_1)\cdots \frac{d\PP}{d\QQ}(x_n) = e^{\sum_{i=1}^n \ln \frac{d\PP}{d\QQ}(x_i)}.
$$
$\PP^{\otimes\infty}$ will stand for the product measure $\PP\otimes\PP\otimes\cdots$ on $[0,1]^\NN$. Taking now the $x_i$'s to be i.i.d. \rvs with common law $\PP$, it follows from the law of large numbers that we have $\PP^{\otimes \infty}$-almost surely
$$
\frac{\PP^{\otimes n}(X_1=x_1,..., X_n=x_n)}{\QQ^{\otimes n}(X_1=x_1,..., X_n=x_n)} \underset{n,\infty}{\simeq} e^{n\EE_{\PP}\bigl[\ln \frac{d\PP}{d\QQ}\bigr]},
$$
in a sense that should be made more precise. The above estimate roughly means that the support of the probability $\PP^{\otimes n}$ in $[0,1]^n$ has $\QQ^{\otimes n}$-measure of order $e^{-n\EE_{\PP}\bigl[\ln \frac{d\PP}{d\QQ}\bigr]}$, when $n$ is large.

One owes to the statistician Charles Stein a rephrasing of this fact in terms of tests, which should be clear from the above description\footnote{Consult for example section $11.7$ of the book \cite{CoverThomas} of Cover and Thomas for a proof of this lemma.}.

\begin{lem}[Stein]
Let $X_1,..., X_n$ be i.i.d. $[0,1]$-valued \rvs with common law $P$. Consider the hypotheses "$H_0 : P=\PP$", and "$H_1 : P=\QQ$", and suppose we want to test hypothesis $H_0$ against $H_1$. The quality of a decision region $A_n\subset [0,1]^n$ is measured by the errors  $\PP^{\otimes n}(A_n^c)$ and $\QQ^{\otimes n}(A_n)$. Given $\varep>0$, set $\beta_n^\varep = \inf\bigl\{\QQ^{\otimes n}(A_n)\,;\,A_n\subset [0,1]^n,\,\PP^{\otimes n}(A_n^c)<\varep\bigr\}$. Then we have
$$
\underset{\varep,0}{\lim}\;\underset{n\infty}{\overline{\lim}}\;\frac{1}{n}\log \beta_n^\varep = -H(\PP\,;\,\QQ).
$$
\end{lem}

To understand this lemma, imagine you want to test the hypothesis "$H_0 : P=\PP$", with a given (very) small bound on the two errors. Then, the smaller $H(\PP\,;\,\QQ)$ will be, the bigger $n$ will have to be in order to design a test achieving the requirements on errors.

If now $\PP$ and $\QQ$ depend on some 'time' $s$ and $H(\PP_s\,;\,\QQ_s)$ decreases, then you will need more and more data to achieve the test, as time passes. It could be said of a situation where $H(\PP_s\,;\,\QQ_s)$ decreases to $0$ that \textit{the process $\bigl(X_1(s), X_2(s)...\bigr)$ forgets its law} as time increases, as it is more and more difficult to distiguish if it has common distribution $\PP_s$ or $\QQ_s$.

\smallskip

From that point of view, a satisfying $H$-theorem for the $(V,\ff)$-processes would be completely different from theorem \ref{ThmDivergence}. Given a $\mcC^1$ path $\gamma$ in $\MMM$ and a (proper) time $s$ (of $\gamma$), define
$$
h_s^\gamma(\be_0\,;\,\bfg) = \frac{q(\dot\gamma_s,\bfg^0)f\bigl(\be_0\,;\,(\gamma_s,\bfg)\bigr)}{\int_{\OO_{\gamma_s}\MMM} q(\dot\gamma_s,\widehat\bfg^0)f\bigl(\be_0\,;\,(\gamma_s,\widehat\bfg)\bigr)\,\Vol_{\gamma_s}(d\widehat\bfg)}.
$$
Let $\VV$ be any (small) spacelike hypersurface \st $\gamma_s\in\VV$ and $\bigl(\dot\gamma_s\bigr)^\perp = T_{\gamma_s}\VV$. The probability $h_s^\gamma\bigl(\be_0\,;\,(\gamma_s,\bfg)\bigr)\Vol_{\gamma_s}(d\bfg)$ is the conditional law of the \rv $\be_{H_{\VV}}{\bf 1}_{H_\VV<\infty}$, given that $m_{H_\VV}=\gamma_s$; we have seen in section \ref{SubSubSectionOnePartFunction} that this conditional probability does not depend on $\VV$ but only on $\gamma_s$ and $\dot\gamma_s$. Given two initial conditions $\be_0, \be_0'$ of the $(V,\ff)$-diffusion, define the time-dependent relative entropy associated to the path $\gamma$ as
$$
H^\gamma_{\be_0,\be_0'}(s) := \int_{\OO_{m_s}\MMM} h_s^\gamma(\be_0\,;\,\bfg) \ln\frac{h_s^\gamma(\be_0\,;\,\bfg)}{h_s^\gamma(\be_0'\,;\,\bfg)}\Vol_{\gamma_s}(d\bfg).
$$
A satisfying $H$-theorem would take the form of the following conjecture.

\begin{conj}
The $\HH\MMM$-valued $(V,\ff)$-diffusion process $(\gamma_.,\gamma'_.) = \bigl\{(m_s,\bfg^0_s)\bigr\}_{s\geq 0}\in\HH\MMM$ almost surely forgets its law as time increases: the relative entropy $H^\gamma_{\be_0,\be_0'}(s)$  decreases $\PP_{\be_0}$-almost surely to $0$, for any $\be_0, \be_0'\in\OO\MMM$.
\end{conj}

It would also be interesting to see if the following holds.

\begin{conj}
A freely falling observer has more and more difficulties in distinguishing $\PP_{\be_0}$ from $\PP_{\be_0'}$.
\end{conj}

We shall adress these questions in a near future. As a last comment, let us notice that theorem \ref{ThmDivergence} can be given an information theoretic flavour. Consider indeed that each open set of spacetime initially has a quantity of "information" equal to its volune, and that this "information" travels with the flow of the vector field $X$. Then, theorem \ref{ThmDivergence} means that the quantity of information that can be found in a fixed open set decreases as the flow-time increases. Yet, this interpretation is far from being as clear as the above two conjectures.

\section{Comments}

It is now time to forget the details of the proofs and summarize the main ideas and results exposed above.

A general class of relativistic diffusions was first presented in the article \cite{ROUPUnifying2}. Although the authors only consider dynamics in Minkowski spacetime, their class of processes is essentially the same as the above class of $(V,\ff)$-diffusions. It is characterized by the existence at each time of a rest frame with the property that the moving object has, in addition to a deterministic accelearation, a Brownian acceleration in any spacelike direction of the rest frame, when computed using the time of the rest frame. Yet, the authors' analysis of the situation rests entirely on a transport equation; an approach similar in spirit to the semi-group analysis of Markov processes, as opposed to the pathwise study of the process. We propose in this article a simple and direct construction of relativistic diffusions on any Lorentzian manifold as flows of stochastic differential equations. This construction necessitates to build the diffusions in the orthonormal frame bundle of $(\MMM,q)$, as was done by Malliavin or Elworthy for Brownian motion in a Riemannian manifold, and by Franchi and Le Jan in the Lorentzian framework. This change of framework is worth being made. Not only are we able to recover directly many of the results established so far, but this pathwise approach presents several other advantages over the analytical method used up to now.

\smallskip

First, it provides a direct (co-ordinate free) description of the dynamics in $\OO\MMM$ which is given as the fundamental mathematical object of the model. Simple hypotheses can be given (section \ref{SubsubsectionSubDiff}) to construct a diffusion in the more familiar phase space $\HH\MMM$ from the diffusion on $\OO\MMM$. An interesting outcome of this approach is the clear new definition of the one-particle distribution function that can be given using the pathwise behaviour of the $(V,\ff)$-process (proposition/definition \ref{DefnOnePartDistr}). The fundamental equation it satisfies (theorem \ref{LStarOnePartFunction}) provides a dynamical justification of the approach used up to now, and sheds some light on the study of the Poisson and Martin boundaries of the $(0,\be_.)$-process (section \ref{SubSubSectionHarmonicFunctions}). Last, but not least, the formalism of vector fields enables us to give in section \ref{BouquetFinal} a concise and clear proof of a general $H$-theorem.

\smallskip

Although these results are encouraging, it would be desirable to discuss the adequacy of the models provided by $(V,\ff)$-diffusions to situations of physical interest. To paraphrase what was written in  the introduction of the seminal article \cite{ROUP0}, the models provided by $(V,\ff)$-diffusions\footnote{By the R.O.U.P. in this article.} should not be considered as accurate models of motion of a "colloidal particle immersed in a real (relativistic) medium". Rather, they should be considered as toy models designed to provide a framework for the study of the main characteristics of the diffusion phenomenon. In this direction, it would certainly be useful to develop an approach to the relativistic Boltzmann equation in terms of hydrodynamic limit of a system of interacting particles\footnote{Consult the article \cite{IsraelBoltzmann} for a discussion of Boltzmann equation in a relativistic framework. The article \cite{AnderssonComer} of Andersson and Comer is also a valuable source of information on relativistic fluid dynamics from a macorscopic point of view.}. Propagation of chaos results could justify the use of $(V,\ff)$-diffusions as models of diffusion dynamics. Other dynamics, as the one introduced by L. Markus in the article \cite{Markus}, might happen to be of physical relevance. Note also the interest that the possibility to define $(V,\ff)$-diffusions in non-isotropic media might have.

Nevertheless, one can consider as one of the merits of our approach the fact that it provides new questions. 
A few of them have been written under the form of conjectures in sections \ref{SubSubSectionHarmonicFunctions} and \ref{BouquetFinal}; we would like to put forwards two other problems concerning the $(0,\be_.)$-process, as we think this is a fundamental object.

\paragraph{Lifetime.} 
The question of explosion of a general $(V,\ff)$-diffusion may appear irrelevant from a physical point of view, after reading the above comments. Yet, the study of this problem for the 'geometric' $(0,\be_.)$-diffusion might happen to be extremely fruitful in its possible links with the existence of singularities of the spacetime itself. Indeed, all the studies made so far, in Minkowski, Robertson-Walker, Schwarzschild and G\"odel spacetimes\footnote{In \cite{BailleulPoisson}, \cite{BailleulRaugi}, \cite{FranchiLeJan} and \cite{Franchi}.} tend to reinforce the feeling that the $\MMM$-part of the $(0,\be_.)$-diffusion eventually behaves like a lightlike geodesic\footnote{The almost-sure convergence of the $\RR^{1,d}$-part of the $(0,\be_.)$-process to a random point of the causal boundary of $\RR^{1,d}$, proved in \cite{BailleulRaugi}, is the clearest proof of this fact.}. So its seems natural to ask the following question.

\begin{OpenPb}
Is null geodesic incompleteness equivalent to explosion of the $(0,\be_.)$-process with positive probability?
\end{OpenPb}
This link between geometry and probability would provide a new approach to the existence of singularities on Lorentzian manifolds. It would be interesting for instance to see wether the hypotheses of Penrose's theorem\footnote{Consult the original article \cite{HawkingPenrose} of Hawking and Penrose or the books \cite{Oneill} for instance.} are relevant from a probabilistic point of view or not. One of its potential benefits is that the explosion problem has an analytical counterpart which is a \textit{linear} problem. Explosion is equivalent to any of the following two conditions\footnote{See the article \cite{Grigoryan} of A. Grigor'yan.}.

\begin{itemize}
   \item[\textit{\textbf{1.}}] Let $\la>0$. There exists a non-null bounded smooth function $f$ \st $(L-\la)f = 0$.
   \item[\textit{\textbf{2.}}] Let $T>0$. There exists a non-null solution to the Dirichlet problem $\partial_t f = Lf$, on $[0,T]\times\MMM$, with initial condition $0$.
\end{itemize}
The use of the one-particle distribution function ${\overset{\leftarrow}{f^\la}}(\be_0\,;\,\cdot)$ of the $\bigl(0,\overset{\leftarrow}{\be}_.\bigr)$-process \textit{killed at constant rate $\la$} will certainly help to see if condition \textit{\textbf{1}} holds. To begin with, it would be interesting to find an example of a geodesically timelike complete Lorentzian manifold whose $(0,\be_.)$-diffusion explodes. No such manifold has been found yet.

\paragraph{A probabilistic interpretation of Einstein tensor?} 
In so far as the local geometry of spacetime can be recovered from the pathwise behaviour of the $(0,\be_.)$-process\footnote{As was noticed in the example following definition \ref{DefnDiff}.}, it is tempting to ask if one can ultimately give a probabilistic interpretation of Einstein tensor determining matter in terms of $(0,\be_.)$-diffusion. This question brings us far from the present day knowledge... We hope it will have some day a positive answer.

\bigskip

\noindent \textbf{Aknowledgements.} I would like to thank Jacques Franchi for his numerous comments on an early version of the manuscript; they led to a clearer exposition of the results exposed here.

\bibliographystyle{unsrt}
\bibliography{biblio}

\end{document}